\documentclass[a4paper,11pt,final]{article}

\usepackage{graphicx,psfrag}
\usepackage{amssymb,amsmath,amsthm,mathrsfs}
\usepackage{enumitem}
\usepackage{a4wide}
\usepackage{hyperref}
\usepackage{array}

\numberwithin{equation}{section}
\allowdisplaybreaks[4]

\newtheorem{proposition}{Proposition}[section]
\newtheorem{theorem}[proposition]{Theorem}
\newtheorem{lemma}[proposition]{Lemma}

\newtheorem{definition}[proposition]{Definition}

\newenvironment{proofof}[1]{\smallskip\noindent{\textbf{Proof~of~#1.}}%
  \hspace{1pt}}{\hspace{-5pt}{\nobreak\quad\nobreak\hfill\nobreak%
    $\square$\vspace{2pt}\par}\smallskip\goodbreak}

\newcommand{\Id}{\mathinner{\mathrm{Id}}}
\newcommand{\pint}[1]{\mathaccent23{#1}}

\newcommand{\C}[1]{\mathbf{C}^{#1}}
\newcommand{\Cc}[1]{\mathbf{C}_c^{#1}}

\newcommand{\BV}{\mathbf{BV}}

\renewcommand{\L}[1]{{\mathbf{L}^#1}}
\newcommand{\Lloc}[1]{{\mathbf{L}_{loc}^{#1}}}

\newcommand{\W}[2]{{\mathbf{W}^{#1,#2}}}

\newcommand{\modulo}[1]{{\left|#1\right|}}
\newcommand{\norma}[1]{{\left\|#1\right\|}}
\newcommand{\caratt}[1]{{\chi_{\strut#1}}}
\newcommand{\reali}{{\mathbb{R}}}
\newcommand{\naturali}{{\mathbb{N}}}

\renewcommand{\epsilon}{\varepsilon}
\renewcommand{\phi}{\varphi}
\renewcommand{\theta}{\vartheta}

\newcommand{\tv}{\mathinner{\rm TV}}
\newcommand{\spt}{\mathop{\rm spt}}
\newcommand{\sgn}{\mathop{\rm sgn}}

\renewcommand{\d}[1]{\mathinner{\mathrm{d}{#1}}}

\newcommand{\vv}{$\boldsymbol{(v)}$}

\newcommand{\Caption}[1]{\caption{{\small#1}}}
\renewcommand{\div}{\nabla \cdot}

\begin{document}

\title{Hyperbolic Predators vs.~Parabolic Preys}

\author{Rinaldo M.~Colombo$^1$ \and Elena Rossi$^2$}

\footnotetext[1]{Unit\`a INdAM, Universit\`a di Brescia,
Italy.}

\footnotetext[2]{Dipartimento di Matematica ed Applicazioni, Universit\`a di Milano-Bicocca, 
Italy.}

\date{Deember 16, 2013}

\maketitle

\begin{abstract}

  \noindent We present a nonlinear predator--prey system consisting of
  a nonlocal conservation law for predators coupled with a parabolic
  equation for preys. The drift term in the predators' equation is a
  nonlocal function of the prey density, so that the movement of
  predators can be directed towards region with high prey
  density. Moreover, Lotka--Volterra type right hand sides describe
  the feeding. A theorem ensuring existence, uniqueness, continuous
  dependence of weak solutions and various stability estimates is
  proved, in any space dimension. Numerical integrations show a few
  qualitative features of the solutions.

  \medskip

  \noindent\textit{2000~Mathematics Subject Classification:} 35L65,
  35M30, 92D25

  \medskip

\end{abstract}

\section{Introduction}
\label{sec:Intro}

Consider the following predator--prey model
\begin{equation}
  \label{eq:Model}
  \left\{
    \begin{array}{l}
      \partial_t u
      +
      \div \left(
        u \, v (w)
      \right)
      =
      (\alpha\,w - \beta) \, u
      \\
      \partial_t w - \mu\, \Delta w
      =
      (\gamma - \delta\, u) \, w
    \end{array}
  \right.
\end{equation}
where $u = u (t,x)$, respectively $w = w (t,x)$, is the predator,
respectively prey, density at time $t \in \reali^+$ and position $x
\in \reali^n$. Preys diffuse according to a parabolic equation, since
$\mu >0$. Here, $\gamma$ is the prey birth rate and $\delta$ the prey
mortality due to the predators. The predator density evolves according
to a hyperbolic balance law, where the coefficient $\alpha$ in the
source term accounts for the increase in the predator density due to
feeding on preys, while $\beta$ is the predator mortality rate. The
flow $u\, v (w)$ accounts for the preferred predators' direction. The
velocity $v$ is in general a \emph{nonlocal} and nonlinear function of
the prey density. A typical choice can be
\begin{equation}
  \label{eq:w}
  v (w)
  =
  \kappa\,
  \frac{\nabla (w*\eta)}{\sqrt{1+\norma{\nabla (w*\eta)}^2}}
\end{equation}
meaning that predators move towards regions of higher concentrations
of preys. Indeed, when $\eta$ is a positive smooth mollifier with
$\int_{\reali^n} \eta \d{x}=1$, the space convolution $\left(w (t) *
  \eta\right) (x)$ has the meaning of an average of the prey density
at time $t$ around position $x$.  The denominator
$\sqrt{1+\norma{\nabla (w*\eta)}^2}$ is merely a smooth normalization
factor, so that the positive parameter $\kappa$ is the maximal
predator speed.

Two key features of the model~\eqref{eq:Model} are the
following. First, while preys diffuse in all directions due to the
Laplacian in the $w$ equation, predators in~\eqref{eq:Model} have a
directed movement, for instance drifting towards regions where the
prey density is higher. This allows, for instance, to describe
predators chasing preys. Second, predators have a well defined
horizon. Indeed, the radius of the support of $\eta$ in~\eqref{eq:w}
defines how far predators can \emph{``feel''} the presence of preys
and, hence, the direction in which they move.

Aim of this paper is to study the class of models~(\ref{eq:Model})
under suitable assumptions on $v$. We prove below existence,
uniqueness, continuous dependence from the initial datum and various
stability estimates for the solutions to~\eqref{eq:Model}. Here,
solutions are found in the space $\L1 \cap \L\infty \cap \BV$ for the
predators and in $\L1 \cap \L\infty$ for the preys. Thus, solutions
are here understood in the distributional sense, see
definitions~\ref{def:Main}, \ref{def:para}
and~\ref{def:Hyp}. Moreover, all analytical results hold in any space
dimension, the explicit dependence of the constants entering the
estimates being duly reported in the proofs below. Besides,
qualitative properties of solutions are shown by means of numerical
integrations.

\medskip

With reference to possible biological applications, the words
\emph{prey} and \emph{predator} should be here understood in their
widest sense. The diffusion in the second equation may well describe
the evolution of a chemical substance or also of temperature. Indeed,
setting for instance $\delta = 0$, the second equation decouples from
the first and the first one fits
into~\cite[Formula~(0.1)]{FriedmanConservation}, see
also~\cite{ColomboGaravelloLecureux, ColomboHertyMercier}. In this
connection, we recall that the interest in nonlocal hyperbolic models
is increasing in several fields.

Various multiD models devoted to crowd dynamics are considered
in~\cite[Section~4]{ColomboHertyMercier} and
in~\cite{ColomboGaravelloLecureux} in the case of a single population,
in~\cite{ColomboLecureuxPerDafermos} for several populations. In these
works, solutions are understood in the weak sense of Kru\v zkov,
see~\cite{Kruzkov}, and well posedness is proved in any space
dimension.

Nonlocal models for aggregation and swarming are presented
in~\cite{Fetecau2011, FetecauEftimie2010}, where the existence of
smooth or Lipschitz continuous solutions is proved in 1D and in 2D,
the $n$ dimensional case being considered in~\cite{FetecauNonlin}. Due
to the biological motivation, in these papers only one population is
considered.

In structured population biology, the use of nonlocal models based on
conservation laws is very common, also in a measure valued setting,
see for example~\cite{AcklehEtAl, CarrilloEtAl, DiekmannGetto} and the
references therein.

On the other hand, the use of purely parabolic equations in
predator--prey models with spatial distributions is rather classical,
see for instance~\cite[Section~1.2]{MurrayII}. With respect to these
models, the use of a first order differential operator in the predator
density allows to describe the directed movement of predators and
ensures that they have a finite propagation speed. Indeed, if the
initial distribution of predators has compact support, then the region
they occupy grows with finite speed and remains compact for all times,
as proved below.

\medskip

As analytical tools, in this paper we consider separately the
equations
\begin{equation}
  \label{eq:due}
  \partial_t u + \div \left(c (t,x) \, u\right) = b (t,x) \, u
  \qquad \mbox{ and } \qquad
  \partial_t w - \mu \, \Delta w = a (t,x) \, w  \,.
\end{equation}
For the former, we exploit the classical results by Kru\v
zkov~\cite{Kruzkov} and the more recent stability estimates proved
in~\cite{ColomboMercierRosini, MercierV2}. The literature on the
latter equation in~\eqref{eq:due} is vast, however our considering it
in $\L1 \cap \L\infty$ on all $\reali^n$ seems to be somewhat
unconventional, hence we provide detailed proofs of the necessary
estimates. The two equations~\eqref{eq:due} are here studied following
exactly the same template and analogous results are obtained. Once the
necessary estimates for the solutions to~\eqref{eq:due} are proved, a
fixed point argument allows to prove the well posedness
of~\eqref{eq:Model} and several stability estimates.

\medskip

The next section presents the analytical results, first the main
theorem and then the propositions at its basis. Section~\ref{sec:NI}
is devoted to sample numerical integrations of~\eqref{eq:Model}. All
technical details are deferred to the final Section~\ref{sec:TD}.

\section{Analytical Results}
\label{sec:AR}

This paragraph is devoted to the well posedness theorem that
constitutes the main result of this paper. All proofs are deferred to
Section~\ref{sec:TD}.

Our first step is the rigorous definition of solution
to~\eqref{eq:Model}.

\begin{definition}
  \label{def:Main}
  Let $T>0$ be fixed. A solution to the system~\eqref{eq:Model} on
  $[0,T]$ is a pair $(u,w) \in \C0 ([0,T]; \L1 (\reali^n;\reali^2))$
  such that
  \begin{itemize}

  \item setting $a (t,x) = \gamma-\delta\,u (t,x)$, $w$ is a weak
    solution to $\partial_t w - \mu \, \Delta w = a \, w$;

  \item setting $b (t,x) = \alpha \, w (t,x) - \beta$ and $c(t,x) =
    \left(v \left(w (t)\right)\right) (x)$, $u$ is a weak solution to
    $\partial_t u + \div \left( u \, c \right) = b \, u$.

  \end{itemize}
\end{definition}

\noindent The extension to the case of the Cauchy problem is
immediate. Below, in Definition~\ref{def:para}, respectively in
Definition~\ref{def:Hyp}, we state and use different definitions of
solutions to the parabolic equation $\partial_t w - \mu \, \Delta w =
a \, w$, respectively to the hyperbolic equation $\partial_t u + \div
\left( u \, c \right) = b \, u$, and prove their equivalence in
Lemma~\ref{lem:paraSol}, respectively in Lemma~\ref{lem:hypSol}.

Throughout, we work in the spaces
\begin{equation}
  \label{eq:X}
  \begin{array}{rcl@{\;}c@{\;}ll}
    \mathcal{X}
    & = &
    (\L1 \cap \L\infty \cap \BV) \left(\reali^n; \reali\right)
    & \times &
    (\L1 \cap \L\infty) \left(\reali^n; \reali\right)
    & \mbox{and}
    \\[3pt]
    \mathcal{X}^+
    & = &
    (\L1 \cap \L\infty \cap \BV) \left(\reali^n; \reali^+\right)
    & \times &
    (\L1 \cap \L\infty) \left(\reali^n; \reali^+\right)
  \end{array}
\end{equation}
with the norm
\begin{equation}
  \label{eq:normaX}
  \norma{(u,w)}_{\mathcal{X}}
  =
  \norma{u}_{\L1 (\reali^n;\reali)}
  +
  \norma{w}_{\L1 (\reali^n;\reali)} \,.
\end{equation}
System~\eqref{eq:Model} is defined by a few real parameters and by the
map $v$, which is assumed to satisfy the following condition:
\begin{description}
\item[\vv] $v \colon (\L1 \cap \L\infty) (\reali^n; \reali) \to (\C2
  \cap \W1\infty (\reali^n; \reali^n))$ admits a constant $K$ and an
  increasing map $C \in \Lloc\infty (\reali^+; \reali^+)$ such that,
  for all $w,w_1, w_2 \in (\L1\cap\L\infty) (\reali^n; \reali)$,
  \begin{eqnarray*}
    \norma{v (w)}_{\L\infty (\reali^n; \reali^n)}
    & \leq &
    K \, \norma{w}_{\L1 (\reali^n; \reali)}
    \\
    \norma{\nabla v (w)}_{\L\infty (\reali^n;\reali^{n\times n})}
    & \leq &
    K \, \norma{w}_{\L\infty(\reali^n; \reali)}
    \\
    \norma{v (w_1) - v (w_2)}_{\L\infty (\reali^n; \reali^n)}
    & \leq &
    K \, \norma{w_1 - w_2}_{\L1(\reali^n; \reali)}
    \\
    \norma{\nabla \left(\div v (w)\right)}_{\L1 (\reali^n; \reali^n)}
    & \leq &
    C \left(\norma{w}_{\L1(\reali^n; \reali)}\right)
    \, \norma{w}_{\L1(\reali^n; \reali)}
    \\
    \norma{\div\left(v (w_1) - v (w_2)\right)}_{\L1 (\reali^n; \reali)}
    & \leq &
    C \left(\norma{w_2}_{\L\infty(\reali^n; \reali)}\right)
    \,
    \norma{w_1 - w_2}_{\L1(\reali^n; \reali)}.
  \end{eqnarray*}
\end{description}

\noindent Above, the bound on the $\L\infty$ norm of $v (w)$ by means
of the $\L1$ norm of $w$ is typical of a nonlocal, e.g.~convolution,
operator. Indeed, Lemma~\ref{lem:v} below ensures that under
reasonable regularity conditions on the kernel $\eta$, the operator
$v$ in~\eqref{eq:w} satisfies~\vv.

Relying solely on~\vv, we state the main result of this paper.

\begin{theorem}
  \label{thm:main}
  Fix $\alpha, \beta, \gamma, \delta \geq 0$ and $\mu > 0$. Assume
  that $v$ satisfies~\vv. Then, there exists a map
  \begin{displaymath}
    \mathcal{R} \colon \reali^+ \times \mathcal{X}^+ \to \mathcal{X}^+
  \end{displaymath}
  with the following properties:
  \begin{enumerate}[leftmargin=*]
    \setlength{\itemsep}{1pt} \setlength{\parskip}{1pt}
  \item \label{it:thm:semi} \textbf{$\boldsymbol{\mathcal{R}}$ is a
      semigroup:} $\mathcal{R}_0 = \Id$ and $\mathcal{R}_{t_2} \circ
    \mathcal{R}_{t_1} = \mathcal{R}_{t_1+t_2}$ for all $t_1,t_2 \in
    \reali^+$.
  \item \label{it:thm:solution} \textbf{$\boldsymbol{\mathcal{R}}$
      solves~\eqref{eq:Model}:} for all $(u_o,w_o) \in \mathcal{X}^+$,
    the map $t \to \mathcal{R}_t (u_o,w_o)$ solves the Cauchy Problem
    \begin{displaymath}
      \left\{
        \begin{array}{l}
          \partial_t u
          +
          \div \left(
            u \, v (w)
          \right)
          =
          (\alpha\,w - \beta) \, u
          \\
          \partial_t w - \mu\, \Delta w
          =
          (\gamma - \delta\, u) \, w
          \\
          u (0,x) = u_o (x)
          \\
          w (0,x) = w_o (x)
        \end{array}
      \right.
    \end{displaymath}
    in the sense of Definition~\ref{def:Main}. In particular, for all
    $(u_o, w_o) \in \mathcal{X}^+$ the map $t \to \mathcal{R}_t (u_o,
    w_o)$ is continuous in time.

  \item \label{it:thm:Lipschitz} \textbf{Local Lipschitz continuity in
      the initial datum:} for all $r > 0$ and for all $t \in
    \reali^+$, there exist a positive $\mathcal{L} (t,r)$ such that
    for all $(u_1,w_1), \, (u_2,w_2) \in \mathcal{X}^+$ with
    \begin{displaymath}
      \norma{u_i}_{\L\infty (\reali^n;\reali)} \leq r
      \,,\qquad
      \tv (u_i) \leq r
      \,,\qquad
      \norma{w_i}_{\L\infty (\reali^n;\reali)} \leq r
      \,,\qquad
      \norma{w_i}_{\L1 (\reali^n;\reali)} \leq r
    \end{displaymath}
    for $i=1,2$, the following estimate holds:
    \begin{displaymath}
      \norma{
        \mathcal{R}_t (u_1,w_1) -
        \mathcal{R}_t (u_2,w_2)}_{\mathcal{X}}
      \leq
      \mathcal{L} (t,r) \, \norma{(u_1,w_1) - (u_2,w_2)}_{\mathcal{X}} \,.
    \end{displaymath}

  \item \label{it:thm:growth} \textbf{Growth estimates:} for all
    $(u_o, w_o) \in \mathcal{X}^+$ and for all $t \in \reali^+$,
    denote $(u,w) (t) = \mathcal{R}_t (u_o,w_o)$. Then,
    \begin{eqnarray*}
      \norma{u (t)}_{\L1 (\reali^n; \reali)}
      & \leq &
      \norma{u_o}_{\L1 (\reali^n; \reali)} \,
      \exp\left(\alpha \frac{e^{\gamma \, t} - 1}{\gamma}
        \norma{w_o}_{\L\infty (\reali^n; \reali)}\right)
      \\
      \norma{u (t)}_{\L\infty (\reali^n; \reali)}
      & \leq &
      \norma{u_o}_{\L\infty (\reali^n; \reali)} \,
      \exp\left((\alpha + K) \frac{e^{\gamma \, t} - 1}{\gamma}
        \norma{w_o}_{\L\infty (\reali^n; \reali)}\right)
      \\
      \norma{w (t)}_{\L1 (\reali^n; \reali)}
      & \leq &
      \norma{w_o}_{\L1 (\reali^n; \reali)} \, e^{\gamma \, t}
      \\
      \norma{w (t)}_{\L\infty (\reali^n; \reali)}
      & \leq &
      \norma{w_o}_{\L\infty (\reali^n; \reali)} \, e^{\gamma \, t} \,.
    \end{eqnarray*}

  \item \label{thm:support} \textbf{Propagation speed:} if $(u_o, w_o)
    \in \mathcal{X}^+$ is such that $\spt (u_o) \subseteq B (0,
    \rho_o)$, then for all $t \in \reali^+$,
    \begin{displaymath}
      \spt \left(u (t)\right) \subseteq  B\left(0, \rho (t) \right)
      \qquad \mbox{ where } \qquad
      \rho (t) = \rho_o + K \, t \, e^{\gamma t} \norma{w_o}_{\L1 (\reali^n;\reali)}\,.
    \end{displaymath}
  \end{enumerate}
\end{theorem}

\noindent An explicit estimate of the Lipschitz constant $\mathcal{L}
(t,r)$ is provided at~\eqref{eq:Lip}.

Theorem~\ref{thm:main} is proved through careful estimates on the
parabolic problem
\begin{equation}
  \label{eq:CauchyP}
  \left\{
    \begin{array}{l}
      \partial_t w - \mu \, \Delta w = a (t,x) \, w
      \\
      w (t_o,x) = w_o (x)
    \end{array}
  \right.
\end{equation}
and, separately, on the balance law
\begin{equation}
  \label{eq:CauchyH}
  \left\{
    \begin{array}{l}
      \partial_t u + \div \left(c (t,x) \, u\right) = b (t,x) \, u
      \\
      u (t_o,x) = u_o (x) \,.
    \end{array}
  \right.
\end{equation}
Our approaches to both the evolution equations~\eqref{eq:CauchyP}
and~\eqref{eq:CauchyH} are identical. We recall below the key
definitions, prove the basic well posedness results and provide
rigorous stability estimates, always referring to the spaces
in~\eqref{eq:X} and with reference to the $\L1$ norm.

To improve the readability of the statements below, we denote by
$\mathcal{O} (t)$ an increasing smooth function of time $t$, depending
on the space dimension $n$ and on various norms of the coefficients
$a, \mu$ in~\eqref{eq:CauchyP} and $b,c$ in~\eqref{eq:CauchyH}. All
proofs are deferred to Section~\ref{sec:TD}, where explicit estimates
for all constants are provided.

Throughout, we fix $t_o, T \in \reali^+$, with $T > t_o$, and denote
\begin{equation}
  \label{eq:IJ}
  I = [t_o,T]
  \quad \mbox{ and } \quad
  J = \left\{(t_1,t_2) \in I^2 \colon t_1 < t_2\right\} .
\end{equation}

For completeness, we recall the following notions from the theory of
parabolic equations. They are similar to various results in the wide
literature on parabolic problems, see for instance~\cite{Amann,
  LunardiBook, QuittnerSouplet}, but here we are dealing with $\L1$
solutions on the whole space.

Inspired by~\cite[Section~48.3]{QuittnerSouplet}, we give the
following definition, where we used the notation~\eqref{eq:IJ}.

\begin{definition}
  \label{def:para}
  Let $a \in \L\infty (I \times \reali^n; \reali)$ and $w_o \in \L1
  (\reali^n; \reali)$.  A \emph{weak solution} to~\eqref{eq:CauchyP}
  is a function $w \in \C0 (I; \L1 (\reali^n; \reali))$ such that for
  all test functions $\phi \in \Cc\infty (\pint{I} \times \reali^n;
  \reali)$
  \begin{equation}
    \label{eq:Psol}
    \int_{t_o}^T \int_{\reali^n}
    (w \, \partial_t \phi + \mu \, w \, \Delta \phi + a \, w \, \phi) \,
    \, \d{x} \, \d{t}
    = 0
  \end{equation}
  and $w(t_o,x) = w_o (x) $.
\end{definition}

The following lemma is similar to various results in the literature,
see for instance~\cite[Section~48.3]{QuittnerSouplet}, and is here
recalled for completeness. The heat kernel is denoted by $H_\mu (t,x)
= (4 \, \pi \, \mu \, t)^{-n/2} \; \exp \left(-\norma{x}^2 \middle/(4
  \, \mu \, t)\right)$, where $t > 0$, $x \in \reali^n$ and $\mu >0$
is fixed.

\begin{lemma}
  \label{lem:paraSol}
  Let $a \in \L\infty (I \times \reali^n; \reali)$. Assume that $w_o
  \in \L1 (\reali^n; \reali)$. Then,
  \begin{enumerate}[leftmargin=*]
  \item any function $w$ satisfying
    \begin{equation}
      \label{eq:solNucleo}
      w (t,x) = \left( H_\mu (t-t_o) * w_o \right) (x)
      +
      \int_{t_o}^t
      \left( H_\mu (t-\tau) * \left(a (\tau) \, w (\tau)\right) \right)
      (x) \, \d\tau
    \end{equation}
    solves~\eqref{eq:CauchyP} in the sense of
    Definition~\ref{def:para};
  \item any solution to~\eqref{eq:CauchyP} in the sense of
    Definition~\ref{def:para} satisfies~\eqref{eq:solNucleo}.
  \end{enumerate}
\end{lemma}

\noindent The well posedness of~\eqref{eq:CauchyP} is now proved.

\begin{proposition}
  \label{prop:para}
  Let $a \in \L\infty (I \times \reali^n; \reali)$. Then, the Cauchy
  problem~\eqref{eq:CauchyP} generates a map $\mathcal{P} \colon J
  \times \L1 (\reali^n; \reali) \to \L1 (\reali^n; \reali)$ with the
  following properties:
  \begin{enumerate}[leftmargin=*]
    \setlength{\itemsep}{1pt} \setlength{\parskip}{1pt}
  \item \textbf{$\boldsymbol{\mathcal{P}}$ is a Process:}
    $\mathcal{P}_{t,t} = \Id$ for all $t \in I$ and
    $\mathcal{P}_{t_2,t_3} \circ \mathcal{P}_{t_1,t_2} =
    \mathcal{P}_{t_1,t_3}$ for all $t_1, t_2, t_3 \in I$, with $t_1
    \leq t_2 \leq t_3$.

  \item \textbf{$\boldsymbol{\mathcal{P}}$
      solves~(\eqref{eq:CauchyP}):} for all $w_o \in \L1 (\reali^n;
    \reali)$, the function $t \to \mathcal{P}_{t_o,t} w_o$ solves the
    Cauchy problem~\eqref{eq:CauchyP} in the sense of
    Definition~\ref{def:para}.

  \item \textbf{Regularity in $\boldsymbol{w_o}$:} for all $(t_o,t)
    \in J$, the map $\mathcal{P}_{t_o,t} \colon \L1 (\reali^n; \reali)
    \to \L1 (\reali^n; \reali)$ is linear and continuous, with
    \begin{displaymath}
      \norma{\mathcal{P}_{t_o,t} w_o}_{\L1 (\reali^n; \reali)}
      \leq
      \mathcal{O} (t) \, \norma{w_o}_{\L1 (\reali^n; \reali)} \, .
    \end{displaymath}

  \item \textbf{$\L\infty$ estimate:} for all $w_o \in (\L1\cap
    \L\infty) (\reali^n; \reali)$, for all $(t_o,t) \in J$,
    \begin{displaymath}
      \norma{\mathcal{P}_{t_o,t}w_o}_{\L\infty (\reali^n; \reali)}
      \leq
      \mathcal{O} (t) \, \norma{w_o}_{\L\infty (\reali^n; \reali)} \, .
    \end{displaymath}

  \item \textbf{Stability in $\boldsymbol{a}$:} let $a_1, a_2 \in
    \L\infty (I \times \reali^n; \reali)$ with $a_1-a_2 \in \L1 (I
    \times \reali^n; \reali)$ and call $\mathcal{P}^1, \mathcal{P}^2$
    the corresponding processes.  Then, for all $(t_o,t) \in J$ and
    for all $w_o \in (\L1\cap \L\infty) (\reali^n; \reali)$,
    \begin{displaymath}
      \norma{
        \mathcal{P}^1_{t_o,t} w_o - \mathcal{P}^2_{t_o,t} w_o}_{\L1 (\reali^n; \reali)}
      \leq
      \mathcal{O} (t) \,
      \norma{w_o}_{\L\infty (\reali^n; \reali)} \,
      \norma{a_1-a_2}_{\L1 ([t_o,t] \times \reali^n;\reali)} \,.
    \end{displaymath}

  \item \textbf{Positivity:} if $w_o \in (\L1 \cap \L\infty)
    (\reali^n; \reali)$ and $w_o \geq 0$, then $\mathcal{P}_{t_o,t} \,
    w_o \geq 0$ for all $(t_o,t) \in J$.

  \item \textbf{Regularity in $\boldsymbol{(t,x)}$:} if $w_o \in (\L1
    \cap \C1) (\reali^n; \reali)$, then $(t,x) \to
    (\mathcal{P}_{t_o,t} w_o) (x) \in \C1 (I \times \reali^n;
    \reali)$.

  \item \textbf{Regularity in time:} for all $w_o \in \L1 (\reali^n;
    \reali)$, the map $t \to \mathcal{P}_{t_o,t} w_o$ is in $\C0
    \left(I; \L1 (\reali^n; \reali)\right)$, moreover for every
    $\theta \in \left]0,1\right[$ and for all $\tau, \, t_1, \, t_2
    \in I$ with $t_2 \geq t_1 \geq \tau > t_o$
    \begin{displaymath}
      \norma{
        \mathcal{P}_{t_o,t_2}w_o - \mathcal{P}_{t_o,t_1}w_o}_{\L1 (\reali^n;\reali)}
      \leq
      \norma{w_o}_{\L1 (\reali^n; \reali)}\,
      \left[
        \frac{n}{\tau-t_o} + \mathcal{O}(t_2)
      \right] \,
      \modulo{t_2-t_1}^\theta.
    \end{displaymath}

  \item \textbf{$\W11$ estimate:} for all $w_o \in \L1 (\reali^n;
    \reali)$, for all $(t_o,t) \in J$,
    \begin{displaymath}
      \norma{\nabla (\mathcal{P}_{t_o,t}w_o)}_{\L1 (\reali^n; \reali^n)}
      \leq
      \frac{\mathcal{O} (t)}{\sqrt{t-t_o}} \, \norma{w_o}_{\L1 (\reali^n; \reali)}.
    \end{displaymath}
  \end{enumerate}
\end{proposition}

\noindent We now follow the same template used in the preceding
proposition and lemma, but referring to the hyperbolic
problem~\eqref{eq:CauchyH}. Similarly
to~\cite[Section~4.3]{DafermosBook}
and~\cite[Section~3.5]{SerreBooks}, we give the following definition,
where we used the notation~\eqref{eq:IJ}.

\begin{definition}
  \label{def:Hyp}
  Let $b \in \L\infty (I \times \reali^n; \reali)$, $c \in \L\infty (I
  \times \reali^n; \reali^n)$ and $u_o \in ( \L1 \cap \L\infty )
  (\reali^n; \reali)$. A \emph{weak solution} to~\eqref{eq:CauchyH} is
  a function $u \in \C0 (I; \L1 (\reali^n; \reali))$ such that for all
  test functions $\phi \in \Cc\infty (\pint{I} \times
  \reali^n;\reali)$
  \begin{equation}
    \label{eq:Hsol}
    \int_{t_o}^T \int_{\reali^n}
    (u \, \partial_t \phi + u \, c \cdot \nabla \phi + b \, u \, \phi) \,
    \d{x} \, \d{t} = 0
  \end{equation}
  and $u (t_o,x) = u_o (x)$.
\end{definition}

The following Lemma is analogous to Lemma~\ref{lem:paraSol}, with the
usual integral formula~\eqref{eq:solNucleo} replaced by integration
along characteristics, see~\eqref{eq:solCara}.

\begin{lemma}
  \label{lem:hypSol}
  Let $c$ be such that $ c \in (\C0\cap \L\infty) (I \times \reali^n;
  \reali^n)$, $c (t) \in \C1 (\reali^n; \reali^n) \quad \forall t \in
  I$, $\nabla c \in \L\infty (I \times \reali^n; \reali^{n\times n})$.
  Assume that $b \in \L\infty (I \times \reali^n; \reali)$ and $u_o
  \in (\L1 \cap \L\infty) (\reali^n; \reali)$. Then
  \begin{enumerate}[leftmargin=*]
  \item the function $u$ defined by
    \begin{equation}
      \label{eq:solCara}
      u (t,x)
      =
      u_o (X (t_o;t,x)) \,
      \exp\left(
        \int_{t_o}^t \left(
          b(\tau,X (\tau;t,x))
          -
          \div c \left(\tau,X (\tau;t,x)\right)
        \right)
        \d\tau
      \right),
    \end{equation}
    where the map $t \mapsto X (t;t_o,x_o)$ solves the Cauchy Problem
    \begin{equation}
      \label{eq:2}
      \left\{
        \begin{array}{l}
          \dot X = c (t, X)
          \\
          X (t_o) = x_o \,,
        \end{array}
      \right.
    \end{equation}
    is a Kru\v zkov solution to~\eqref{eq:CauchyH}, i.e. for all $k
    \in \reali$ and for all $\phi \in \Cc\infty (\pint{I} \times
    \reali^n; \reali^+)$,
    \begin{equation}
      \label{eq:Ksol}
      \int_{t_o}^T \int_{\reali^n}
      \left[
        (u - k) (\partial_t \phi
        +
        c \cdot \nabla \phi)
        +
        (b \,u - k \, \div c) \, \phi
      \right]
      \sgn (u - k) \, \d{x} \, \d{t} \geq 0
    \end{equation}
    and $u (t_o,x) = u_o (x)$, hence, $u$ solves~\eqref{eq:CauchyH} in
    the sense of Definition~\ref{def:Hyp};
  \item any solution to~\eqref{eq:CauchyH} in the sense of
    Definition~\ref{def:Hyp} coincides with $u$ as defined
    in~\eqref{eq:solCara}.
  \end{enumerate}
\end{lemma}

\noindent We now prove the well posedness of~\eqref{eq:CauchyH}.

\begin{proposition}
  \label{prop:hyper}
  We pose the assumptions:
  \begin{description}
  \item[(b)] $b \in (\C1 \cap \L\infty)(I \times \reali^n; \reali)$;
    $\nabla b \in \L1 (I \times \reali^n; \reali^n)$;
  \item[(c)] $c \in (\C2 \cap \L\infty) (I \times \reali^n;
    \reali^n)$; $\nabla c \in \L\infty (I \times \reali^n;
    \reali^{n\times n})$; $\nabla (\div c) \in \L1(I \times
    \reali^n;\reali^n)$.
  \end{description}
  Then, the Cauchy Problem~\eqref{eq:CauchyH} generates a map
  $\mathcal{H} \colon J \times (\L1 \cap \L\infty \cap \BV) (\reali^n;
  \reali) \to (\L1 \cap \L\infty \cap \BV) (\reali^n; \reali)$ with
  the following properties:
  \begin{enumerate}[leftmargin=*]
    \setlength{\itemsep}{1pt} \setlength{\parskip}{1pt}
  \item \textbf{$\boldsymbol{\mathcal{H}}$ is a process:}
    $\mathcal{H}_{t,t} = \Id$ for all $t \in I$ and
    $\mathcal{H}_{t_2,t_3} \circ \mathcal{H}_{t_1,t_2} =
    \mathcal{H}_{t_1,t_3}$ for all $t_1, t_2, t_3 \in I$, with $t_1
    \leq t_2 \leq t_3$.

  \item \textbf{$\boldsymbol{\mathcal{H}}$ solves~(\ref{eq:CauchyH}):}
    for all $u_o \in (\L1 \cap \L\infty \cap \BV) (\reali^n; \reali)$,
    the function $t \to \mathcal{H}_{t_o,t} u_o$ solves the Cauchy
    problem~\eqref{eq:CauchyH} in the sense of
    Definition~\ref{def:Hyp}.

  \item \textbf{Regularity in $\boldsymbol{u_o}$:} for all $(t_o,t)
    \in J$ the map $\mathcal{H}_{t_o,t} \colon (\L1 \cap \L\infty \cap
    \BV) (\reali^n; \reali) \to (\L1 \cap \L\infty \cap \BV)
    (\reali^n; \reali)$ is linear and continuous, with
    \begin{displaymath}
      \norma{\mathcal{H}_{t_o,t} u_o}_{\L1 (\reali^n; \reali)}
      \leq
      \mathcal{O} (t) \, \norma{u_o}_{\L1 (\reali^n; \reali)}
      \,.
    \end{displaymath}

  \item \textbf{$\L\infty$ estimate:} for all $u_o \in (\L1 \cap
    \L\infty \cap \BV) (\reali^n; \reali)$, for all $(t_o,t) \in J$,
    \begin{displaymath}
      \norma{\mathcal{H}_{t_o,t} u_o}_{\L\infty (\reali^n; \reali)}
      \leq
      \mathcal{O} (t) \, \norma{u_o}_{\L\infty (\reali^n; \reali)}
      \,.
    \end{displaymath}

  \item \textbf{Stability in $\boldsymbol{b,c}$:} if $b_1, b_2$
    satisfy~\textbf{(b)} with $b_1 -b_2 \in \L1 (I\times \reali^n;
    \reali)$ and $c_1, c_2$ satisfy~\textbf{(c)} with $\div (c_1 -c_2)
    \in \L1 (I \times \reali^n; \reali)$, call $\mathcal{H}^1,
    \mathcal{H}^2$ the corresponding processes. Then, for all $(t_o,t)
    \in J$ and for all $u_o \in (\L1 \cap \L\infty \cap \BV)
    (\reali^n; \reali)$,
    \begin{eqnarray*}
      & &
      \norma{
        \mathcal{H}^1_{t_o,t} u_o - \mathcal{H}^2_{t_o,t} u_o}_{\L1 (\reali^n; \reali)}
      \\
      & \leq &
      \mathcal{O} (t)
      \left( \norma{u_o}_{\L\infty (\reali^n; \reali)} + \tv (u_o) \right)
      \norma{c_1-c_2}_{\L1 ([t_o,t];\L\infty (\reali^n; \reali^n))}
      \\
      & &
      +
      \mathcal{O} (t) \norma{u_o}_{\L\infty (\reali^n; \reali)}
      \left(
        \norma{b_1-b_2}_{\L1 ([t_o,t]\times\reali^n; \reali)}
        +
        \norma{\div (c_1-c_2)}_{\L1 ([t_o,t]\times\reali^n; \reali)}
      \right) .
    \end{eqnarray*}

  \item \textbf{Positivity:} if $u_o \in (\L1\cap\L\infty\cap \BV)
    (\reali^n; \reali)$ and $u_o \geq 0$, then $\mathcal{H}_{t_o, t}
    \, u_o \geq 0$ for all $(t_o,t) \in J$.

  \item \textbf{Total variation bound:} if $u_o \in
    (\L1\cap\L\infty\cap \BV) (\reali^n; \reali)$, then, for all
    $(t_o,t) \in J$,
    \begin{displaymath}
      \tv\left(\mathcal{H}_{t_o,t} u_o \right)
      \leq
      \mathcal{O} (t)
      \left(\norma{u_o}_{\L\infty (\reali^n; \reali)} + \tv (u_o) \right) \,.
    \end{displaymath}

  \item \textbf{Regularity in time:} for all $u_o \in (\L1 \cap
    \L\infty \cap \BV) (\reali^n; \reali)$, the map $t \to
    \mathcal{H}_{t_o,t} u_o$ is in $\C{0,1} \left(I; \L1(\reali^n;
      \reali)\right)$, moreover for all $t_1, t_2 \in J$,
    \begin{displaymath}
      \norma{\mathcal{H}_{t_o,t_2} u_o - \mathcal{H}_{t_o,t_1} u_o}_{\L1 (\reali^n; \reali)}
      \leq
      \mathcal{O} (t_2) \!
      \left(
        \norma{u_o}_{\L1 (\reali^n; \reali)}
        +
        \norma{u_o}_{\L\infty (\reali^n; \reali)}
        +
        \tv (u_o)
      \right)\!
      \modulo{t_2-t_1} .
    \end{displaymath}

  \item \textbf{Finite propagation speed:} let $(t_o, t) \in J$ and
    $u_o \in (\L1 \cap \L\infty \cap \BV) (\reali^n; \reali)$ have
    compact support $\spt u_o$. Then, also $\spt
    \mathcal{H}_{t_o,t}u_o$ is compact.

  \end{enumerate}
\end{proposition}

\section{Numerical Integrations}
\label{sec:NI}

To illustrate some qualitative properties of the solutions
to~\eqref{eq:Model}, we present the result of a few numerical
integrations.

To integrate both equations we use the operator splitting algorithm to
combine the differential operators and the source terms. The balance
law is integrated by means of the Lax--Friedrichs scheme with
dimensional splitting~\cite[Section~19.5]{LeVeque}, while its source
term is solved using a second order Runge--Kutta method. For the
parabolic equation, we use the forward finite differences algorithm
and the usual Euler forward explicit method on the source term.  We
leave the proof of the convergence of this algorithm to the
forthcoming work~\cite{RossiSchleper}.

Remark that the numerical integration of~\eqref{eq:Model} requires a
convolution integral to be computed at each time step. This puts a
constraint on the space mesh, which should be sufficiently small with
respect to the radius of the support of the convolution kernel to
allow a good approximation of this integral.

Here, we focus on the two-dimensional case, that is $n = 2$, and use
the vector field $v$ in~\eqref{eq:w} with the compactly supported
mollifier
\begin{equation}
  \label{eq:eta}
  \eta (x)
  =
  \hat\eta \,
  \left(\ell^2 - \norma{x}^2\right)^3 \caratt{B (0,\ell)} (x)
  \quad \mbox{ with } \hat\eta \in \reali^+  \mbox{ such that }\quad
  \int_{\reali^2} \eta (x) \d{x} = 1 \,.
\end{equation}

The analytical theory developed above is referred to the Cauchy
problem on the whole space $\reali^2$. In both examples below, the
numerical domain of integration is the rectangle $[-1,1] \times
[-2,2]$. The necessary boundary conditions are different in the two
cases and are specified below. The time step $(\Delta t)_P$ for the
parabolic equation and the one $(\Delta t)_H$ for the hyperbolic part
are chosen so that $(\Delta t)_P$ is of the order of $\left((\Delta
  t)_H\right)^2$. The time step for the hyperbolic equation complies
with the usual CFL condition.

Below, we constrain both unknown functions $u$ and $w$ to remain equal
to the initial datum all along the boundary, which is acceptable in
the first equation since no wave in the solution to the balance law
ever hits the numerical boundary. Concerning the second equation, the
choice of these boundary conditions amounts to assume that the
displayed solution is part of a solution defined on all $\reali^2$
that gives a constant inflow into the computational domain.

\subsection{Predators Chasing Preys}

We present a situation in which the effect of the first order
transport term in the predator equation is clearly visible, as well as
the well known Lotka--Volterra type effect in which a species
apparently almost disappears and then its density rises again.

\begin{figure}[h!]
  \centering
  \includegraphics[width=0.33\textwidth,trim = 40 0 20
  0]{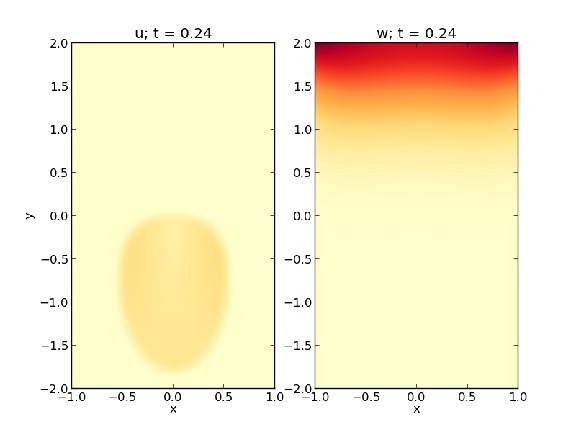}\hfil%
  \includegraphics[width=0.33\textwidth,trim = 40 0 20
  0]{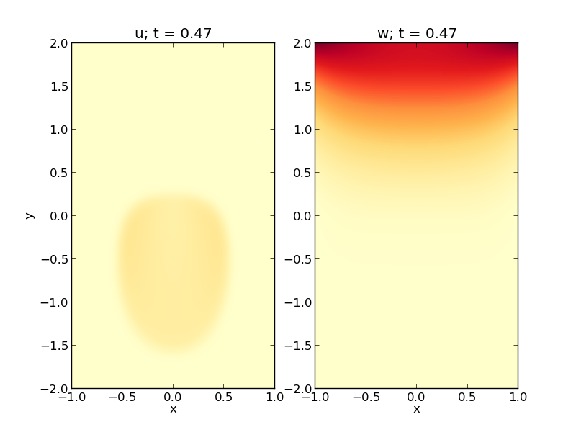}\hfil%
  \includegraphics[width=0.33\textwidth,trim = 40 0 20
  0]{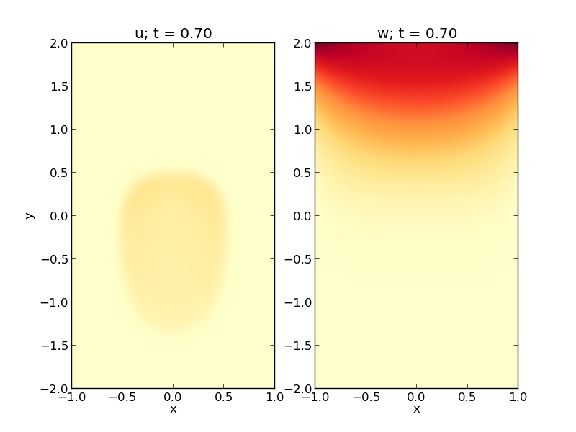}\\
  \includegraphics[width=0.33\textwidth,trim = 40 0 20
  0]{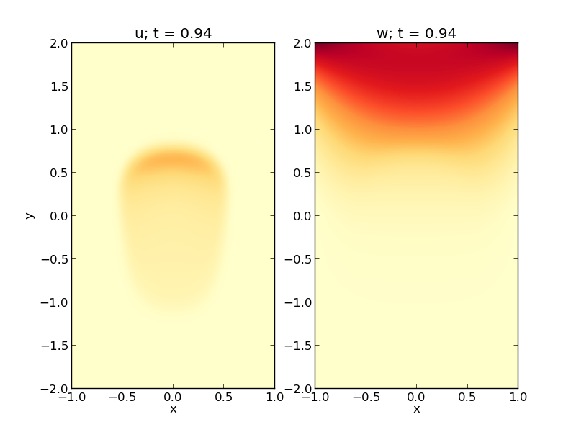}\hfil%
  \includegraphics[width=0.33\textwidth,trim = 40 0 20
  0]{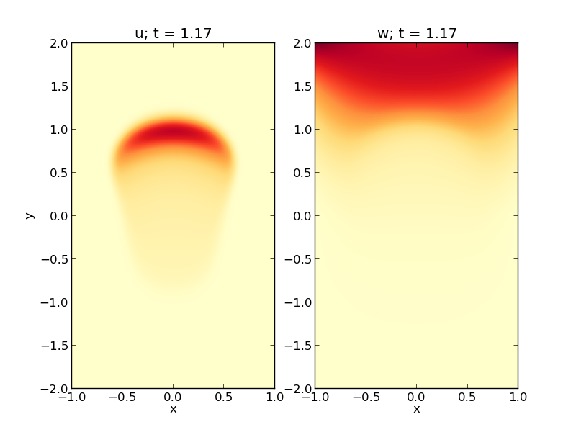}\hfil%
  \includegraphics[width=0.33\textwidth,trim = 40 0 20
  0]{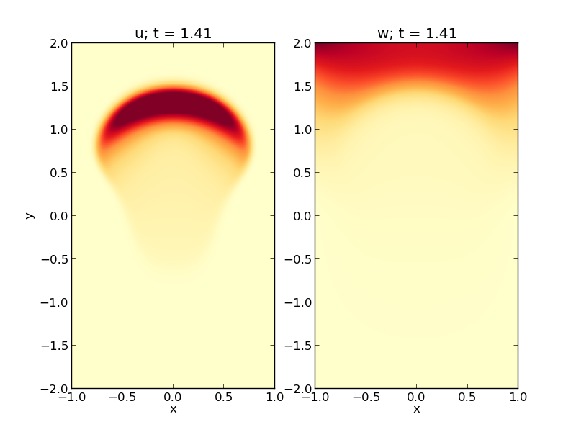}\\
  \Caption{Numerical integration
    of~\eqref{eq:Model}--\eqref{eq:eta}--\eqref{eq:PCPpara} with
    initial datum~\eqref{eq:PCPID} at times $t= 0.24, \, 0.47,\,
    0.70,\, 0.94, \, 1.17, \, 1.41$. In each couple of figures, the
    predator density $u$ is on the left and the prey density $w$ is on
    the right. The colors range as in~\eqref{eq:PCPID} in the interval
    $[0, \, 15]$ for $u$ and $[0, \, 14]$ for $w$. First, predators
    decrease due to lack of nutrients. Thanks to diffusion, preys
    reach the zone where they are \emph{``seen''} by predators.  Then,
    clearly, predators are attracted towards preys and their density
    starts to increase. This solution was obtained with space mesh
    $\Delta x = \Delta y = 0.005$.}
  \label{fig:PCP}
\end{figure}

\begin{figure}[h!]
  \centering
  \includegraphics[width = 0.45\textwidth,trim = 10 0 30 10,
  clip=true]{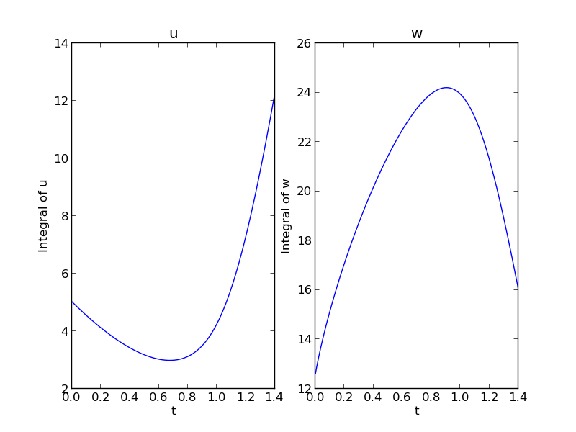}\\
  \Caption{The integrals of $u$, left, and $w$, right, over the
    computational domain versus time; $u$ and $w$ are the numerical
    solutions
    to~\eqref{eq:Model}--\eqref{eq:eta}--\eqref{eq:PCPpara}.}
  \label{fig:PCPcontrol}
\end{figure}

We set $v$ as in~\eqref{eq:w}, $\eta$ as in~\eqref{eq:eta} and
\begin{equation}
  \label{eq:PCPpara}
  \begin{array}{rcl@{\quad\qquad}rcl@{\quad\qquad}rcl}
    \alpha & = & 2
    & \beta & = & 1
    & \kappa & = & 1
    \\
    \gamma & = & 1
    & \delta & = & 2
    & \mu & = & 0.5
  \end{array}
  \qquad \ell \; = \; 0.15
\end{equation}
with initial datum\\
\begin{minipage}[c]{0.4\linewidth}
  \hspace{-0.5pc}
  \begin{tabular}{m{0.85\textwidth}m{0.05\textwidth}@{}}
    \hspace{-0.2pc}
    \includegraphics[width=0.8\textwidth,trim = 40 0 50 0,clip=true]{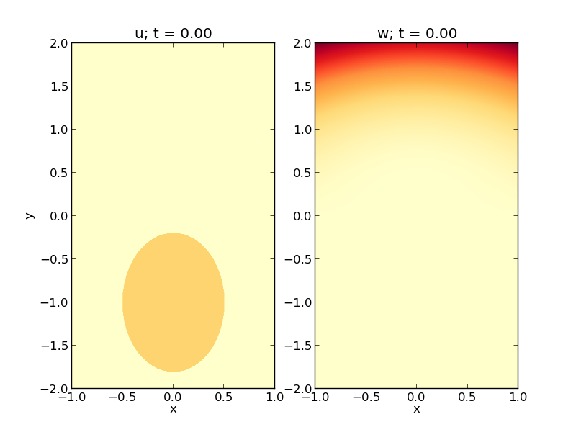} &
    \hspace{-1.2pc}\includegraphics[height=0.16\textheight]{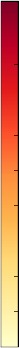}
  \end{tabular}
\end{minipage}%
\begin{minipage}{0.6\linewidth}
  \begin{equation}
    \label{eq:PCPID}
    \begin{array}{@{}l@{}}
      \,u_o (x,y)
      =
      4 \; \caratt{A} \! (x,y)
      \\
      w_o (x,y)
      =
      1.5 \, y \, \max\{ 0, x^2 + y^2 - 0.25\} \,
      \caratt{B} \! (x,y)
      \\
      \mbox{ where}
      \\[3pt]
      A  =  \{(x,y) \in \reali^2 \colon (2 \; x)^2 + \left(1.25 \; (y+1)\right)^2 \leq 1\}
      \\
      B  =   \{(x,y) \in \reali^2 \colon y \geq 0\} \,.
    \end{array}
  \end{equation}
\end{minipage}\\
The result of the numerical integration is in Figure~\ref{fig:PCP}.
At first, preys are outside the horizon of predators. Hence the latter
decrease. Thanks to diffusion, some preys enter the region where
predators feel their presence. This causes predators to move towards
the highest prey density. Therefore, predators immediately increase
and their effect on the prey population is clearly seen, as shown also
by the graph of the integrals of $u$ and $w$ in
Figure~\ref{fig:PCPcontrol}.

We remark that, in the present setting, as time grows, undesired
effects due to the presence of the boundary become relevant.

\subsection{A Dynamic Equilibrium}

\begin{figure}[h!]
  \centering
  \includegraphics[width=0.33\textwidth,trim = 40 0 20 0,
  clip=true]{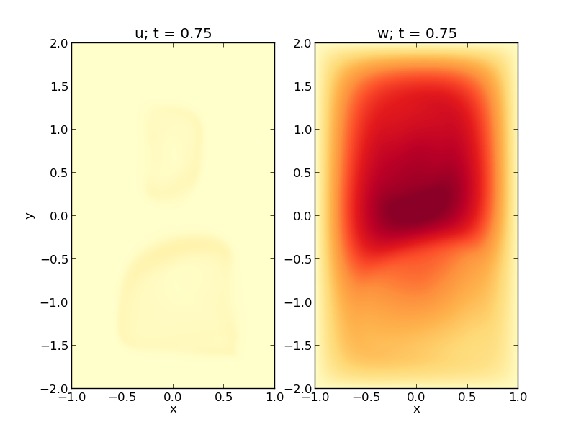}\hfil%
  \includegraphics[width=0.33\textwidth,trim = 40 0 20 0,
  clip=true]{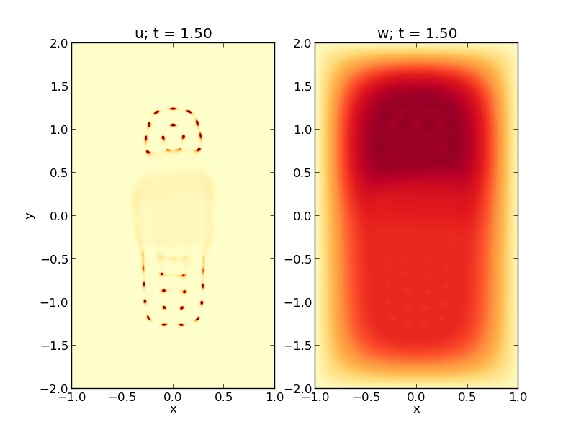}\hfil%
  \includegraphics[width=0.33\textwidth,trim = 40 0 20
  0, clip=true]{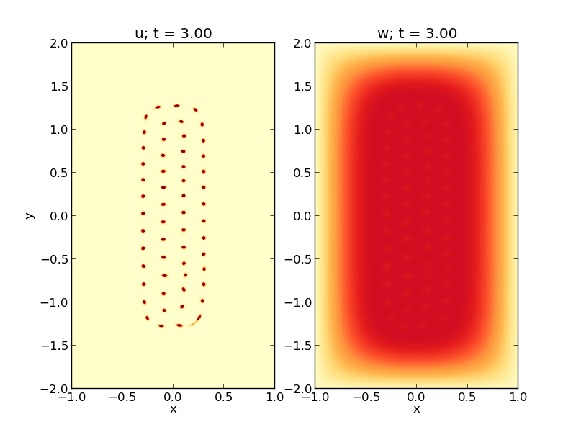}\\
  \includegraphics[width=0.33\textwidth,trim = 40 10 20 0,
  clip=true]{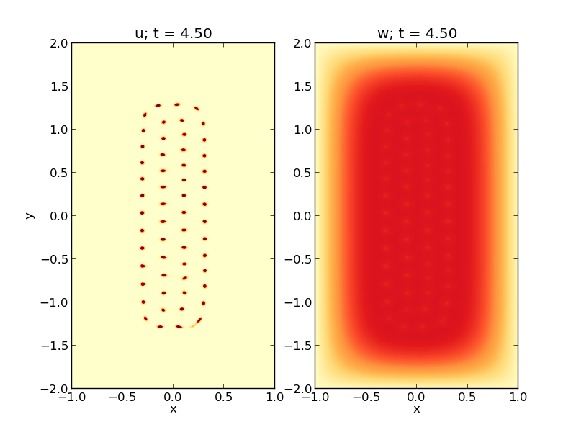}\hfil%
  \includegraphics[width=0.33\textwidth,trim = 40 10 20 0,
  clip=true]{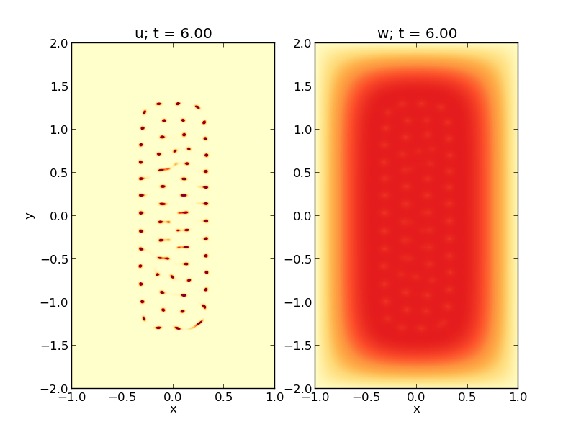}\hfil%
  \includegraphics[width=0.33\textwidth,trim = 40 10 20
  0, clip=true]{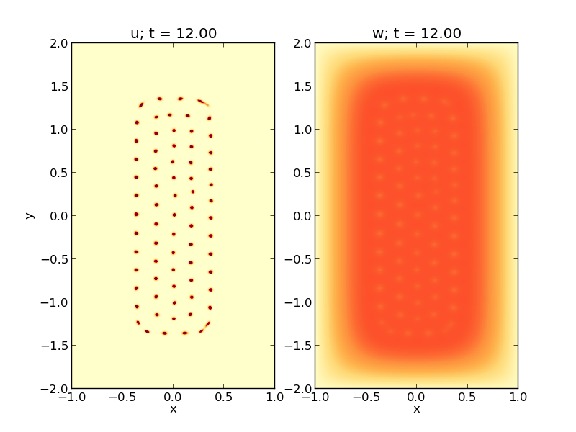}\\
  \Caption{Numerical integration
    of~\eqref{eq:Model}--\eqref{eq:eta}--\eqref{eq:DEpara} with
    initial datum~\eqref{eq:DEID} at times $t= 0.75, \, 1.50, \, 3.00,
    \, 4.50, \, 6.00, \, 12.00$. In each couple of figures, the
    predator density $u$ is on the left and the prey density $w$ is on
    the right. The colors range as in~\eqref{eq:DEID} in the interval
    $[0, \, 0.40]$ for $u$ and $[0.20, \, 0.24]$ for $w$. First,
    predators decrease due to lack of nutrients and move towards the
    central region.  Then, a discrete periodic pattern arises with
    predators focused in small regions regularly distributed along $4$
    columns and, at a later time, along $5$ columns. This solution was
    obtained with space mesh $\Delta x = \Delta y = 0.005$.}
  \label{fig:DE}
\end{figure}

In this case, the numerical solution displays an interesting
asymptotic state in which the diffusion caused by the Laplacian in the
prey equation counterbalances the first order nonlocal differential
operator in the predator equation. The outcome is the onset of a
discrete, quite regular, structure, see Figure~\ref{fig:DE}. We set
$v$ as in~\eqref{eq:w}, $\eta$ as in~\eqref{eq:eta} and
\begin{equation}
  \label{eq:DEpara}
  \begin{array}{rcl@{\quad\qquad}rcl@{\quad\qquad}rcl}
    \alpha & = & 1
    & \beta & = & 0.2
    & \kappa & = & 1
    \\
    \gamma & = & 0.4
    & \delta & = &24
    & \mu & = & 0.5
  \end{array}
  \qquad \ell = 0.25
\end{equation}
with initial datum\\
\begin{minipage}[c]{0.4\linewidth}
  \hspace{-0.5pc}
  \begin{tabular}{m{0.85\textwidth}m{0.05\textwidth}@{}}
    \hspace{-0.2pc}
    \includegraphics[width=0.8\textwidth,trim = 40 0 50 0,clip=true]{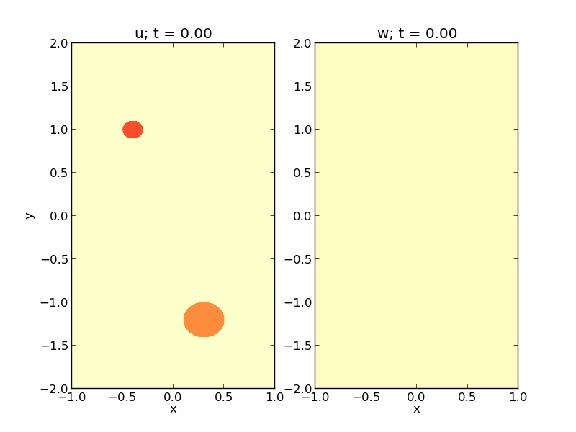} &
    \hspace{-1.2pc}\includegraphics[height=0.16\textheight]{cbarY.jpg}
  \end{tabular}
\end{minipage}%
\begin{minipage}{0.6\linewidth}
  \begin{equation}
    \label{eq:DEID}
    \begin{array}{@{}l@{}}
      \,u_o (x,y)
      =
      0.25 \; \caratt{C} \! (x,y) + 0.2 \; \caratt{D} (x,y)
      \\
      w_o (x,y) = 0.2
      \\[3pt]
      \mbox{ where}
      \\[3pt]
      C \! = \! \{(x,y) \!\in \reali^2 \colon \! (x + 0.4)^2 + (y - 1)^2 \! < 0.01\}
      \\
      D \! = \! \{(x,y)\! \in \reali^2 \colon \! (x - 0.3)^2 + (y + 1.2)^2 \! < 0.04\} \,.
    \end{array}
  \end{equation}
\end{minipage}\\
In this integration, predators first almost disappear, move towards
the central part of the numerical domain and then start to
increase. Slowly, a regular pattern arises. Predators focus in small
regions regularly distributed. These regions display a fairly stable
behavior while passing from being arranged along $4$ to along $5$
columns, see Figure~\ref{fig:DE}, second line.

\begin{figure}[h!]
  \centering
  \includegraphics[width=0.33\textwidth,trim = 40 0 10 0,
  clip=true]{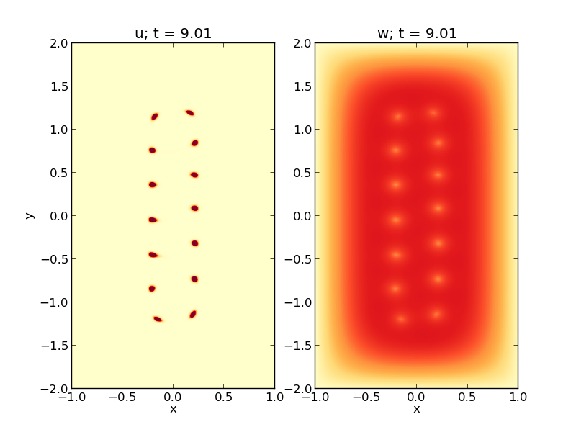}\hfil%
  \includegraphics[width=0.33\textwidth,trim = 40 0 10 0,
  clip=true]{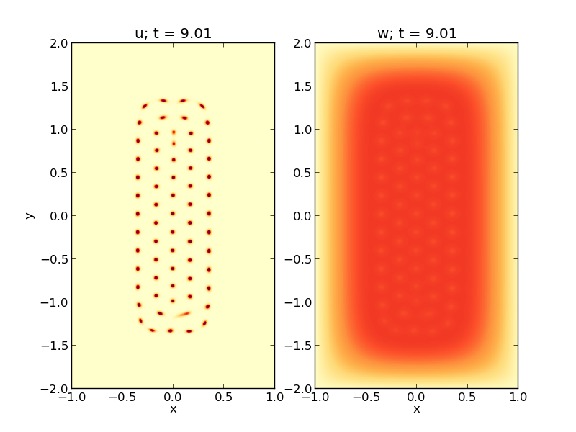}\hfil%
  \includegraphics[width=0.33\textwidth,trim = 40 0 10 0, clip=true]{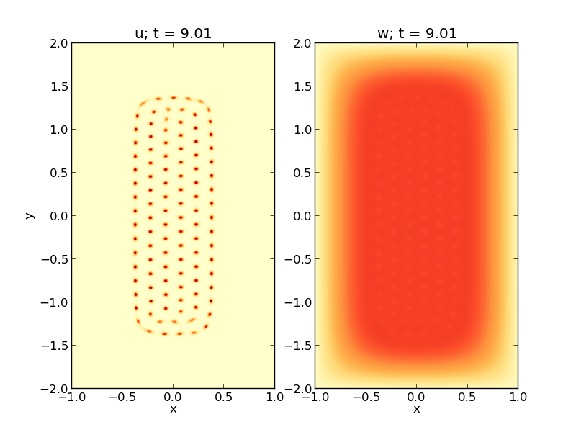}\\
  \Caption{Solutions
    to~\eqref{eq:Model}--\eqref{eq:eta}--\eqref{eq:DEpara} with
    initial datum~\eqref{eq:DEID} computed at time $t = 9.01$ with
    different values of $\ell$, i.e. left $\ell = 0.5$, middle $\ell =
    0.25$ and, right, $\ell = 0.1875$. As $\ell$ decreases, also the
    distance among peaks in the $u$ density decreases and more peaks
    are possible. The color scale is as in Figure~\ref{fig:DE}. These
    solutions were obtained with space mesh $\Delta x = \Delta y =
    0.0075$.}
  \label{fig:DEeyes}
\end{figure}

\begin{figure}[h!]
  \centering
  \includegraphics[width=0.33\textwidth,trim = 40 10 10 0,
  clip=true]{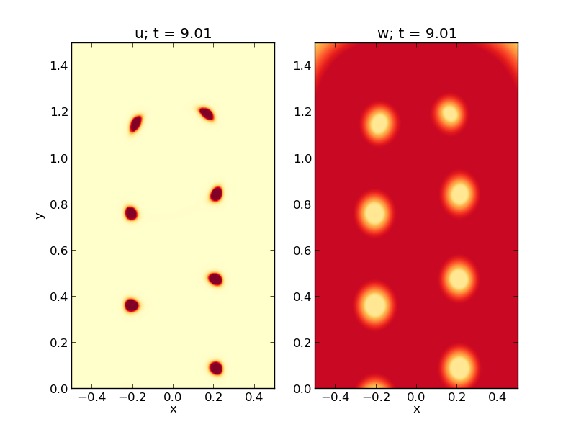}\hfil%
  \includegraphics[width=0.33\textwidth,trim = 40 10 10 0,
  clip=true]{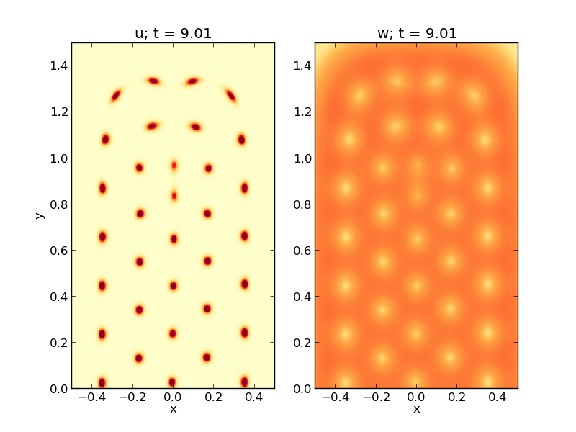}\hfil%
  \includegraphics[width=0.33\textwidth,trim = 40 10 10 0, clip=true]{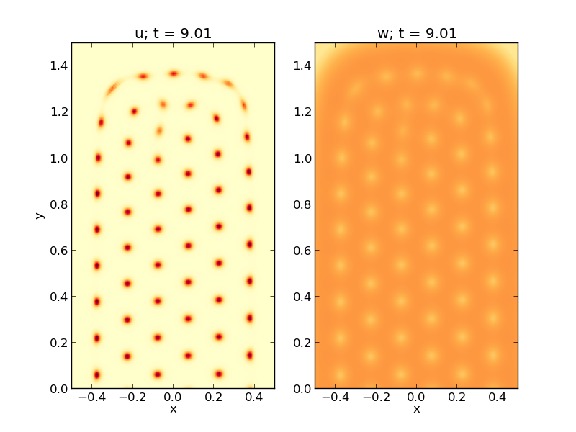}\\
  \Caption{Zoom of Figure~\ref{fig:DEeyes}, restricted to the
    intervals $x \in [-0.5, \, 0.5]$ and $y \in [0, \, 1.5]$. Colors
    are as in~\eqref{eq:DEID}, ranging in $[0, \, 0.4]$ for the $u$
    variable and in $[0.225, \, 0.229]$ for the $w$ variable.}
  \label{fig:DEeyesBis}
\end{figure}

From the analytical point of view, this pattern can be explained as a
dynamic equilibrium between the first order non local operator present
in the predator equation and the Laplacian in the prey equation. Where
predators accumulate, their feeding on preys causes a \emph{``hole''}
in the prey density, see figures~\ref{fig:DEeyes}
and~\ref{fig:DEeyesBis}. As a consequence, the average gradient of the
prey density, which directs the movement of predators, almost vanishes
by symmetry considerations. Hence, predators almost do not move. At
the same time, the diffusion of preys keeps filling the
\emph{``holes''}, thus providing a persistent amount of nutrient to
predators.
Coherently with this explanation, numerical integrations confirm that
the above asymptotic state essentially depends on the size of the
support of $\eta$. Indeed, the mean distance between pairwise nearest
peaks in the density of $u$ is slightly smaller than $\ell$, see
figures~\ref{fig:DEeyes} and~\ref{fig:DEeyesBis}.

\section{Technical Details}
\label{sec:TD}

\subsection{Proofs Related to \texorpdfstring{$\partial_t w - \mu \,
    \Delta u = a (t,x) \, w$}{the Parabolic Equation}}

The following constant will be of use below:
\begin{equation}
  \label{eq:3P}
  J_n = \frac{\Gamma\left((n+1)/2\right)}{\Gamma (n/2)}
\end{equation}
where $\Gamma$ is the Gamma function. Moreover, recall the classical
estimates on the heat kernel:
\begin{eqnarray}
  \label{eq:calore1}
  \norma{H_\mu (t)}_{\L1 (\reali^n; \reali)} & = & 1
  \\
  \label{eq:calore}
  \norma{\nabla H_\mu (t)}_{\L1 (\reali^n; \reali)}
  & = &
  \frac{1}{\sqrt{\mu\,t}} \,
  \frac{\Gamma\left((n+1)/2\right)}{\Gamma (n/2)}
  \; = \;
  \frac{J_n}{\sqrt{\mu\,t}}
  \\
  \nonumber
  \partial_t H_\mu (t,x)
  & = &
  \frac{1}{(4\, \pi \, \mu \, t)^{n/2}} \,
  \frac{\norma{x}^2 - 2\, n \, \mu \, t}{4 \, \mu \, t^2} \,
  \exp\left(-\frac{\norma{x}^2}{4\mu t}\right)
  \\
  \label{eq:StimaHot}
  \norma{\partial_t H_\mu (t)}_{\L1 (\reali^n; \reali)}
  & \leq &
  \frac{n}{t} \,.
\end{eqnarray}

\begin{proofof}{Lemma~\ref{lem:paraSol}}
  \noindent\textbf{1.}
  Since $H_\mu (t) \in \C\infty (\reali^n; \reali)$ and $w$
  satisfies~\eqref{eq:solNucleo}, clearly $w (t) \in \C\infty
  (\reali^n; \reali)$ for all $t > t_o$. If $\phi \in \Cc\infty
  (\pint{I} \times \reali^n; \reali)$ then $\int_{\reali^n} \mu \, w
  \, \Delta \phi \, \d{x} = \int_{\reali^n} \mu \, \Delta w \, \phi \,
  \d{x}$.  Note first that $\mu \, \Delta H_\mu (t) = \partial_t H_\mu
  (t)$ and compute preliminarily
  \begin{eqnarray*}
    \mu \, \Delta w (t,x)
    & = &
    \mu \, \Delta \left( \left( H_\mu
        (t-t_o) * w_o \right) (x) \right) + \int_{t_o}^t \mu \, \Delta
    \left( \left( H_\mu (t-\tau) * \left(a (\tau) \, w (\tau)\right)
      \right) (x) \right) \, \d\tau
    \\
    & = &
    \left( \mu \, \Delta H_\mu (t-t_o) * w_o \right) (x) +
    \int_{t_o}^t \left( \mu \, \Delta H_\mu (t-\tau) * \left(a (\tau)
        \, w (\tau) \right) \right) (x) \, \d\tau
    \\
    & = &
    \left( \partial_t H_\mu (t-t_o) * w_o \right) (x) +
    \int_{t_o}^t \left( \partial_t H_\mu (t-\tau) * \left(a (\tau) \,
        w (\tau) \right) \right) (x) \, \d\tau
    \\
    & = &
    \partial_t \left( H_\mu (t-t_o) * w_o \right) (x) +
    \int_{t_o}^t \partial_t \left( H_\mu (t-\tau) * \left(a (\tau) \,
        w (\tau) \right) \right) (x) \, \d\tau.
  \end{eqnarray*}
  Setting for simplicity $\mathscr{H} (t,x) =
  \displaystyle\int_{t_o}^t\left( H_\mu (t-\tau) * \left(a (\tau) w
      (\tau) \right) \right) (x) \d\tau$, the line above becomes
  \begin{displaymath}
    \mu \, \Delta w (t,x) =
    \partial_t \left( H_\mu (t-t_o) * w_o \right) (x)
    +
    \frac{\d{}}{\d{t}} \mathscr{H} (t,x) - a(t,x) \, w (t,x).
  \end{displaymath}
  We are now able to prove~\eqref{eq:Psol}:
  \begin{eqnarray*}
    & &
    \int_{t_o}^T \int_{\reali^n} (w \, \partial_t \phi + \mu \,
    \Delta w \, \phi + a \, w \, \phi) \, \d{x} \, \d{t}
    \\
    & = &
    \int_{t_o}^T \int_{\reali^n} \Bigl[ \left( H_\mu (t-t_o) *
      w_o \right) (x) \, \partial_t \phi (t,x) + \mathscr {H} (t,x)
    \, \partial_t \phi (t,x)
    \\
    & &
    \qquad +
    \partial_t \left( \left( H_\mu (t-t_o) * w_o \right) (x) \right)\,
    \phi (t,x)
    +
    \left(\frac{\d{}}{\d{t}} \mathscr{H} (t,x)\right) \phi (t,x)
    -
    a(t,x) \, w (t,x) \, \phi (t,x)
    \\
    & &
    \qquad + a(t,x) \, w (t,x) \, \phi (t,x) \Bigr] \, \d{x} \,
    \d{t}
    \\
    & = & \int_{t_o}^T \int_{\reali^n} \frac{\d{}}{\d{t}} \left[
      \left( \left( H_\mu (t-t_o) * w_o \right) (x) + \mathscr {H}
        (t,x) \right) \, \phi (t,x) \right] \, \d{x} \, \d{t}
    \\
    & = & 0.
  \end{eqnarray*}
  It is immediate to verify that if $w$
  satisfies~\eqref{eq:solNucleo}, the initial condition holds.

  \noindent\textbf{2.}
  Let $w$ satisfy~\eqref{eq:solNucleo} and $w_*$ be a weak solution
  to~\eqref{eq:CauchyP}. Then, by the step above, the function $W = w
  - w_*$ is a weak solution to the linear equation~\eqref{eq:CauchyP}
  with zero initial datum.

  Fix $\tau \in \left]t_o,T\right]$, choose any $\phi \in \C1 (I;\Cc1
  (\reali^n;\reali))$ and let $\beta_\epsilon \in \C1 (I;\reali)$ such
  that $\beta_\epsilon (t) = 1$ for $t \in [t_o + \epsilon,\tau -
  \epsilon]$, $\beta_\epsilon' (t) \in [0,2/\epsilon]$ for $t \in
  [t_o,t_o + \epsilon]$, $\beta_\epsilon' (t) \in [-2/\epsilon,0]$ for
  $t \in [\tau - \epsilon,\tau]$ and $\beta_\epsilon (t) = 0$ for $t
  \in [\tau,T]$. Using the definition of weak solution,
  \begin{eqnarray*}
    0
    & = &
    \int_I \int_{\reali^n}
    \left(
      W \, \partial_t (\phi\, \beta_\epsilon)
      +
      \mu W \Delta (\phi \, \beta_\epsilon)
      +
      a \, W \, \phi \, \beta_\epsilon
    \right)
    \d{x} \, \d{t}
    \\
    & = &
    \int_{t_o}^\tau \!\! \int_{\reali^n}
    W \left(
      \partial_t \phi
      +
      \mu \Delta \phi
      +
      a \, \phi \right) \beta_\epsilon \,
    \d{x} \, \d{t}
    +
    \int_{t_o}^{t_o + \epsilon} \!\!\!
    \beta_\epsilon'
    \int_{\reali^n} \!\!\!
    W \, \phi \, \d{x} \, \d{t}
    +
    \int_{\tau-\epsilon}^\tau \!\!\!
    \beta_\epsilon'
    \int_{\reali^n} \!\!\! W \, \phi \, \d{x} \, \d{t}
    \\
    & = &
    \int_{t_o}^\tau
    \beta_\epsilon
    \int_{\reali^n}
    W
    \left(
      \partial_t \phi + \mu \Delta \phi + a \, \phi
    \right) \, \d{x} \, \d{t}
    \\
    & & \quad
    +
    \int_0^1
    \epsilon \, \beta_\epsilon' (\epsilon t + t_o)
    \int_{\reali^n} W (\epsilon t + t_o,x) \, \phi (\epsilon t + t_o,x)\,
    \d{x} \, \d{t}
    \\
    & &
    \quad
    +
    \int_{\tau-\epsilon}^\tau \beta_\epsilon' (t)
    \int_{\reali^n} W (t,x)\, \phi (t,x)\, \d{x} \, \d{t} \,.
  \end{eqnarray*}
  As $\epsilon \to 0$, the first term converges to $\int_{t_o}^\tau
  \int_{\reali^n} W \left( \partial_t \phi + \mu \Delta \phi + a \,
    \phi \right) \, \d{x} \, \d{t}$. By the Dominated Convergence
  Theorem, the second term tends to $0$, since $W (t_o) =
  0$. Concerning the third term, note that
  \begin{eqnarray*}
    & &
    \modulo{
      \int_{\tau-\epsilon}^\tau \beta_\epsilon' (t)
      \int_{\reali^n} W (t,x)\, \phi (t,x)\, \d{x} \, \d{t}
      +
      \int_{\reali^n} W (\tau,x) \, \phi (\tau,x) \, \d{x}
    }
    \\
    & = &
    \modulo{
      \int_{\tau-\epsilon}^\tau
      \beta_\epsilon' (t)
      \int_{\reali^n}
      \left(
        W (t,x)\, \phi (t,x) - W (\tau,x) \, \phi (\tau,x)
      \right)
      \, \d{x} \, \d{t} } \to 0
    \quad \mbox{ as } \quad \epsilon \to 0
  \end{eqnarray*}
  by the continuity of $W$ in time and the smoothness of $\phi$.

  Choose now any $\eta \in \Cc1 (\reali^n;\reali)$ and define $\phi$
  as the backward solution to
  \begin{equation*}
    \left\{
      \begin{array}{l}
        \partial_t \phi + \mu \Delta \phi = - a\, \phi \\
        \phi (\tau) = \eta \,.
      \end{array}
    \right.
  \end{equation*}
  Then,
  \begin{displaymath}
    0
    =
    \int_{t_o}^\tau \int_{\reali^n}
    W
    \left( \partial_t \phi + \mu \Delta \phi + a \, \phi \right)
    \, \d{x} \, \d{t}
    -
    \int_{\reali^n} W (\tau,x) \, \phi (\tau, x) \, \d{x}
    =
    - \int_{\reali^n} W (\tau,x) \, \eta (x) \, \d{x}
  \end{displaymath}
  proving that $W (\tau)$ vanishes identically. By the arbitrariness
  of $\tau$, $w=w_*$.
\end{proofof}

\begin{proofof}{Proposition~\ref{prop:para}}
  For all $t \in I$ denote $A (t) = \exp\left(\int_{t_o}^t \norma{a
      (\tau)}_{\L\infty (\reali^n;\reali)}\d\tau\right)$.

  The proofs of~\textbf{1.}--\textbf{2.} are well known in the
  parabolic literature, see~\cite{Amann, LunardiBook}. By
  Lemma~\ref{lem:paraSol}, recall that the solution $w (t,x) =
  (\mathcal{P}_{t_o,t} w_o) (x)$ to~\eqref{eq:CauchyP}
  satisfies~\eqref{eq:solNucleo}.

  \noindent\textbf{3.}\quad Standard computations, using
  also~\eqref{eq:calore1}, lead to:
  \begin{eqnarray*}
    \norma{w (t)}_{\L1 (\reali^n; \reali)}
    & \leq &
    \int_{\reali^n} \int_{\reali^n}
    H_\mu (t-t_o, x-\xi) \,
    \modulo{w_o (\xi)}
    \d\xi \d x
    \\
    & &
    +
    \int_{\reali^n} \int_{t_o}^t \int_{\reali^n}
    H_\mu (t-\tau, x-\xi) \,
    \modulo{a (\tau,x) \, w (\tau,x)}
    \d\xi \d\tau \d x
    \\
    & \leq &
    \norma{w_o}_{\L1 (\reali^n; \reali)}
    +
    \int_{t_o}^t \norma{a (\tau) \, w (\tau)}_{\L1 (\reali^n; \reali)} \d\tau
    \\
    & \leq &
    \norma{w_o}_{\L1 (\reali^n; \reali)}
    +
    \int_{t_o}^t \norma{a (\tau)}_{\L\infty (\reali^n; \reali)} \,
    \norma{w (\tau)}_{\L1 (\reali^n; \reali)} \d\tau.
  \end{eqnarray*}
  An application of Gronwall Lemma yields the thesis:
  \begin{displaymath}
    \norma{w (t)}_{\L1 (\reali^n; \reali)}
    \leq
    A (t) \, \norma{w_o}_{\L1 (\reali^n; \reali)}.
  \end{displaymath}

  \noindent\textbf{4.}\quad
  By~\eqref{eq:solNucleo} and~\eqref{eq:calore1},
  \begin{eqnarray*}
    \norma{ w (t)}_{\L\infty (\reali^n; \reali)}
    & \leq &
    \norma{w_o}_{\L\infty (\reali^n; \reali)}
    +
    \int_{t_o}^t \norma{a (\tau) \, w (\tau)}_{\L\infty (\reali^n; \reali)} \d\tau
    \\
    & \leq &
    \norma{w_o}_{\L\infty (\reali^n; \reali)}
    +
    \int_{t_o}^t \norma{a (\tau)}_{\L\infty (\reali^n; \reali)}
    \norma{w (\tau)}_{\L\infty (\reali^n; \reali)} \d\tau.
  \end{eqnarray*}
  An application of Gronwall Lemma gives the desired result:
  \begin{displaymath}
    \norma{w (t)}_{\L\infty (\reali^n; \reali)}
    \leq
    A (t) \, \norma{w_o}_{\L\infty (\reali^n; \reali)} \, .
  \end{displaymath}

  \noindent\textbf{5.}\quad Denote $w_i (t) =
  \mathcal{P}^{\phantom{,}i}_{t_o,t}w_o$ and $A_i (t) =
  \exp\left(\int_{t_o}^t \norma{a_i (\tau)}_{\L\infty (\reali^n;
      \reali)}\d\tau\right)$, for $i=1,2$.  Use~\eqref{eq:solNucleo},
  \eqref{eq:calore1} and~4.~above:
  \begin{eqnarray*}
    & &
    \norma{w_1 (t) - w_2 (t)}_{\L1 (\reali^n; \reali)}
    \\
    & \leq &
    \int_{\reali^n} \int_{t_o}^t \int_{\reali^n}
    H_\mu (t-\tau, x-\xi) \,
    \modulo{a_1 (\tau,\xi) \, w_1 (\tau,\xi)
      -
      a_2 (\tau,\xi) \, w_2 (\tau,\xi)} \d\xi \, \d\tau \, \d{x}
    \\
    & \leq &
    \int_{t_o}^t
    \norma{H_\mu (t-\tau)}_{\L1 (\reali^n; \reali)} \,
    \norma{a_1 (\tau) \, w_1 (\tau)-a_2 (\tau) \, w_2 (\tau)}_{\L1 (\reali^n; \reali)} \, \d\tau
    \\
    & = &
    \int_{t_o}^t
    \norma{a_1 (\tau) \, w_1 (\tau)-a_2 (\tau) \, w_2 (\tau)}_{\L1 (\reali^n; \reali)} \, \d\tau
    \\
    & \leq &
    \int_{t_o}^t
    \norma{a_1 (\tau) \, w_1 (\tau)-a_2 (\tau) \, w_1 (\tau)}_{\L1 (\reali^n; \reali)} \, \d\tau
    \! + \!\!
    \int_{t_o}^t \!
    \norma{a_2 (\tau) \, w_1 (\tau)-a_2 (\tau) \, w_2 (\tau)}_{\L1 (\reali^n; \reali)} \, \d\tau
    \\
    & \leq &
    \int_{t_o}^t
    \norma{a_1 (\tau) - a_2 (\tau)}_{\L1 (\reali^n; \reali)} \,
    \norma{w_1 (\tau)}_{\L\infty (\reali^n; \reali)} \d\tau
    \\
    & &
    +
    \int_{t_o}^t
    \norma{a_2 (\tau)}_{\L\infty (\reali^n; \reali)} \,
    \norma{w_1 (\tau)- \, w_2 (\tau)}_{\L1 (\reali^n; \reali)}
    \, \d\tau
    \\
    & \leq &
    \norma{w_o}_{\L\infty (\reali^n; \reali)}
    \int_{t_o}^t
    A_1 (\tau)
    \norma{a_1 (\tau) - a_2 (\tau)}_{\L1 (\reali^n; \reali)} \,
    \d\tau
    \\
    & &
    +
    \int_{t_o}^t
    \norma{a_2 (\tau)}_{\L\infty (\reali^n; \reali)} \,
    \norma{w_2 (\tau)- \, w_1 (\tau)}_{\L1 (\reali^n; \reali)}
    \, \d\tau.
  \end{eqnarray*}
  An application of Gronwall Lemma yields the estimate:
  \begin{displaymath}
    \norma{w_1 (t) - w_2 (t)}_{\L1 (\reali^n; \reali)}
    \leq
    A_1 (t) \, A_2 (t) \,
    \norma{w_o}_{\L\infty (\reali^n; \reali)} \,
    \norma{a_1-a_2}_{\L1 ([t_o,t]\times \reali^n;\reali)}
    \,.
  \end{displaymath}

  \noindent\textbf{6.}\quad Thanks to the $\L\infty$ estimate at~4., we can
  apply~\cite[Chapter~2, Section~4, Theorem~9]{FriedmanBook}.

  \noindent\textbf{7.}\quad As is well known, note that~\eqref{eq:solNucleo}
  immediately ensures that $w$ is of class $\C1$.

  \noindent\textbf{8.}\quad By~\eqref{eq:solNucleo},
  \begin{eqnarray*}
    \norma{w(t_2) - w (t_1)}_{\L1 (\reali^n; \reali)}
    & \leq &
    \norma{
      \left(H_\mu (t_2-t_o) - H_\mu (t_1-t_o)\right)*w_o}_{\L1 (\reali^n; \reali)}
    \\
    & &
    +
    \int_{t_o}^{t_1} \norma{\left(H_\mu (t_2-s) - H_\mu (t_1-s)\right) *
      \left( a (s) w (s)\right)}_{\L1 (\reali^n; \reali)} \, \d{s}
    \\
    & &
    + \int_{t_1}^{t_2}
    \norma{H_\mu (t_2-s) * \left(a (s) w (s)\right)}_{\L1 (\reali^n; \reali)}\d{s}
  \end{eqnarray*}
  and we compute the three terms separately. The first one is the
  usual term of the heat equation, so that using~\eqref{eq:StimaHot},
  \begin{eqnarray*}
    \norma{\left(H_\mu (t_2-t_o) - H_\mu (t_1-t_o)\right)*w_o}_{\L1 (\reali^n; \reali)}
    & \leq &
    \norma{w_o}_{\L1 (\reali^n; \reali)}
    \int_{t_1-t_o}^{t_2-t_o} \norma{\partial_t H_\mu (s)}_{\L1 (\reali^n; \reali)} \, \d{s}
    \\
    & \leq &
    \frac{n}{\tau-t_o} \,
    \norma{w_o}_{\L1 (\reali^n; \reali)} \, \modulo{t_2-t_1}.
  \end{eqnarray*}
  Concerning the second term, use~\eqref{eq:StimaHot}, point~3.~above
  and follow the proof of~\cite[Proposition~4.2.4]{LunardiBook}: for
  every $\theta \in \left]0,1 \right[$ we have
  \begin{eqnarray*}
    & &
    \int_{t_o}^{t_1} \norma{\left(H_\mu (t_2-s) - H_\mu (t_1-s)\right) *
      \left( a (s) w (s)\right)}_{\L1 (\reali^n; \reali)} \, \d{s}
    \\
    & \leq &
    \int_{t_o}^{t_1} \norma{
      \left(\int_{t_1-s}^{t_2-s} \partial_t H_\mu (\sigma) \d\sigma \right) *
      \left( a (s) w (s)\right)}_{\L1 (\reali^n; \reali)} \d{s}
    \\
    & \leq &
    \int_{t_o}^{t_1}
    \norma{\int_{t_1-s}^{t_2-s} \partial_t H_\mu (\sigma) \d\sigma}_{\L1 (\reali^n; \reali)}
    \norma{a (s) w (s)}_{\L1 (\reali^n; \reali)} \d{s}
    \\
    & \leq &
    \int_{t_o}^{t_1}
    \left( \int_{t_1-s}^{t_2-s} \norma{\partial_t H_\mu (\sigma)}_{\L1 (\reali^n; \reali)} \d\sigma \right)
    \norma{a (s)}_{\L\infty (\reali^n; \reali)}
    \norma{w (s)}_{\L1 (\reali^n; \reali)} \d{s}
    \\
    & \leq &
    \int_{t_o}^{t_1}
    \left( \int_{t_1-s}^{t_2-s} \frac{n}{\sigma} \, \d\sigma \right)
    \norma{a (s)}_{\L\infty (\reali^n; \reali)}
    \norma{w_o}_{\L1 (\reali^n; \reali)} \, A (s) \d{s}
    \\
    & \leq &
    \norma{a}_{\L\infty ([t_o,t_1] \times \reali^n; \reali)}
    \norma{w_o}_{\L1 (\reali^n; \reali)} \, A (t_1)
    \int_{t_o}^{t_1} \frac{n}{(t_1-s)^\theta}
    \int_{t_1-s}^{t_2-s} \frac{1}{\sigma^{1-\theta}} \, \d\sigma \d{s}
    \\
    & \leq &
    \norma{a}_{\L\infty ([t_o,t_1] \times \reali^n; \reali)}
    \norma{w_o}_{\L1 (\reali^n; \reali)} \, A (t_1)
    \int_{t_o}^{t_1} \frac{n}{\theta \, (t_1-s)^\theta}
    \left[ (t_2-s)^\theta - (t_1-s)^\theta\right] \, \d{s}
    \\
    & \leq &
    \norma{a}_{\L\infty ([t_o,t_1] \times \reali^n; \reali)}
    \norma{w_o}_{\L1 (\reali^n; \reali)} \, A (t_1)
    \int_{t_o}^{t_1} \frac{n}{\theta \, (t_1-s)^\theta}
    \, \modulo{t_2-t_1}^\theta \, \d{s}
    \\
    & \leq &
    \norma{a}_{\L\infty ([t_o,t_1] \times \reali^n; \reali)}
    \norma{w_o}_{\L1 (\reali^n; \reali)} \, A (t_1) \,
    \frac{n \, (t_1-t_o)^{1-\theta}}{\theta \, (1-\theta)} \,
    \modulo{t_2-t_1}^\theta \,,
  \end{eqnarray*}
  where we used the inequality $(t_2-s)^\theta - (t_1-s)^\theta \leq
  (t_2-t_1)^\theta$.  Concerning the third term,
  use~\eqref{eq:calore1} and point~3.~above
  \begin{eqnarray*}
    \int_{t_1}^{t_2}
    \norma{H_\mu (t_2-s) * \left(a (s) w (s)\right)}_{\L1 (\reali^n; \reali)} \d{s}
    & \!\leq \!&
    \int_{t_1}^{t_2}
    \norma{H_\mu (t_2-s)}_{\L1 (\reali^n; \reali)}
    \norma{a (s) w (s)}_{\L1 (\reali^n; \reali)} \d{s}
    \\
    & \!\leq\! &
    \norma{a}_{\L\infty ([t_1,t_2] \times \reali^n; \reali)} \,
    \norma{w_o}_{\L1 (\reali^n; \reali)} \, A (t_2) \, \modulo{t_2-t_1}.
  \end{eqnarray*}
  Summing up the above expressions, we obtain
  \begin{eqnarray*}
    & &
    \norma{w (t_2) - w (t_1)}_{\L1 (\reali^n; \reali)}
    \\
    & \leq &
    \norma{w_o}_{\L1 (\reali^n; \reali)}\,
    \left[
      \frac{n}{\tau-t_o} +
      \left(\frac{n \, (t_1-t_o)^{1-\theta}}{\theta \, (1-\theta)} +1 \right)\,
      \norma{a}_{\L\infty ([t_0,t_2] \times \reali^n; \reali)} \, A (t_2)
    \right] \,
    \modulo{t_2-t_1}^\theta
  \end{eqnarray*}
  and the H\"older estimate at point~8.~is proved. To prove the
  continuity in time, we are left to check the right continuity in
  $t_o$. To this aim, use~\eqref{eq:calore1} and point~3.~above,
  introduce the variable $\zeta = \xi / (\mu\,t)^{n/2}$ and compute
  \begin{displaymath}
    \begin{array}{@{}c@{\,}l@{}}
      &
      \norma{w (t) - w_o}_{\L1 (\reali^n; \reali)}
      \\
      \leq &
      \displaystyle
      \int_{\reali^n} \! \int_{\reali^n}
      H_\mu (t-t_o, \xi) \modulo{ w_o (x-\xi) - w_o (x)} \d\xi \d{x}
      +
      \int_{t_o}^t
      \norma{H_\mu (t-s) * \left(a (s)\, w (s)\right)}_{\L1 (\reali^n; \reali)}
      \d{s}
      \\
      \leq &
      \displaystyle\int_{\reali^n} \!\! H_1 (1, \zeta) \!\!
      \int_{\reali^n} \!\!
      \modulo{w_o \left(x- (\mu t)^{\frac{n}{2}}\zeta\right) - w_o (x)}
      \d\zeta \d{x}
      +
      A (t) \norma{a}_{\L\infty ([t_o,t]\times\reali^n; \reali)}
      \norma{w_o}_{\L1 (\reali^n;\reali)} \modulo{t-t_o} .
    \end{array}
  \end{displaymath}
  Both terms above vanish as $t \to t_o$, (in the first use the
  Dominated Convergence Theorem), completing the proof of continuity
  in time.

\noindent\textbf{9.}\quad Using~\eqref{eq:calore},
\eqref{eq:solNucleo}, the standard properties of the convolution
product and point~3.~above,
\begin{eqnarray*}
  \!\!\!\!
  \norma{\nabla w (t)}_{\L1 (\reali^n; \reali)}
  \!\!\!\!
  & \leq &
  \!\!\!\!
  \norma{\nabla H_\mu (t-t_o)}_{\L1 (\reali^n; \reali)} \,
  \norma{w_o}_{\L1 (\reali^n; \reali)}
  \\
  \!\!\!\!
  & &
  \!\!\!\!
  +
  \int_{t_o}^t \norma{\nabla H_\mu (t - \tau)}_{\L1 (\reali^n; \reali)}
  \, \norma{a (\tau)}_{\L\infty (\reali^n; \reali)}
  \, \norma{w (\tau)}_{\L1 (\reali^n; \reali)}
  \d\tau
  \\
  \!\!\!\!
  & \leq &
  \!\!\!\!
  \frac{J_n}{\sqrt{\mu \, (t-t_o)}} \,
  \norma{w_o}_{\L1 (\reali^n; \reali)}
  +
  \norma{a}_{\L\infty (I \times\reali^n; \reali)}
  \int_{t_o}^t
  \frac{J_n \, A (\tau)}{\sqrt{\mu (t-\tau)}} \,
  \norma{w_o}_{\L1 (\reali^n; \reali)} \,
  \d\tau
  \\
  \!\!\!\!
  & \leq &
  \!\!\!\!
  \frac{J_n}{\sqrt{\mu \, (t-t_o)}}
  \norma{w_o}_{\L1 (\reali^n; \reali)}
  \left(
    1
    +
    2 \, (t-t_o) \, A (t) \,
    \norma{a}_{\L\infty (I\times \reali^n \reali)}
  \right).
\end{eqnarray*}
\end{proofof}

\subsection{Proofs Related to \texorpdfstring{$\partial_t u + \div
    \left(c (t,x) \, u\right) = b (t,x) \, u$}{the Hyperbolic
    Equation}}

The following constants will be of use below:
\begin{equation}
  \label{eq:3H}
  I_n
  =
  n \int_0^{\pi/2} (\cos \theta)^n \d\theta
  =
  n \,
  \frac{\Gamma\left((n+1)/2\right) \, \Gamma (1/2)}{2 \, \Gamma\left((n+2)/2\right)}
  \,.
\end{equation}

\begin{proofof}{Lemma~\ref{lem:hypSol}}
  We recall~\cite[Section~3]{AmbrosioCrippaUMI} and follow the proof
  of~\cite[Lemma~5.1]{ColomboHertyMercier}.

  \noindent\textbf{1.} Let $u_{o,n} \in \C1 (\reali^n; \reali)$
  approximate $u_o$ in the sense that $\norma{u_{o,n} - u_o}_{\L1
    (\reali^n;\reali)} \to 0$ as $n \to + \infty$.  Call $u_n$ the
  corresponding quantity as given by~\eqref{eq:solCara}. Then,
  $\norma{u_n - u}_{\L\infty(I;\L1 (\reali^n;\reali))} \to 0$ as $n
  \to + \infty$, so that $u \in \L\infty(I;\L1
  (\reali^n;\reali))$. Concerning the continuity in time,
  by~\eqref{eq:solCara} $u_n \in \C0 (I;\L1 (\reali^n;\reali))$ and
  $u$ is the uniform limit of the sequence $u_n$, hence $u \in
  \C0(I;\L1 (\reali^n;\reali))$.

  Using the flow generated by~\eqref{eq:2}, introduce the change of
  variable $y = X (t_o; t, x)$, so that $x = X (t; t_o, y)$. Denote
  its Jacobian by $J (t,y) = \det (\nabla_y X (t; 0, y))$. Then, $J$
  solves
  \begin{displaymath}
    \frac{\d J (t, y)}{\d{t}} = \div c \left(t, X (t; t_o, y)\right) \, J (t,y)
    \quad \mbox{ with } \quad
    J (t_o,y) = 1 \,.
  \end{displaymath}
  Hence
  \begin{displaymath}
    J (t,y)
    =
    \exp \left(\int_{t_o}^t \div c (\tau, X (\tau; t_o,y)) \d\tau \right),
  \end{displaymath}
  so that~\eqref{eq:solCara} can be written as
  \begin{equation}
    \label{eq:CambioVar}
    u (t,x) = \frac{1}{J(t,y)} \, u_o (y) \, \mathcal{B} (t,y)
    \quad \mbox{ where } \quad
    \begin{array}{@{}rcl@{}}
      \mathcal{B} (t,y)
      & = &
      \displaystyle
      \exp\left( \int_{t_o}^t b (\tau,X (\tau;t_o,y)) \d\tau \right)
      \\
      x & = & X (t; t_o, y)
    \end{array}
  \end{equation}
  Let $k \in \reali$ and $\phi \in \Cc\infty (\pint{I} \times
  \reali^n; \reali^+)$. We prove~\eqref{eq:Ksol} for $u$ given as
  in~\eqref{eq:solCara}:
  \begin{eqnarray*}
    \!\!\!\!\!\!
    & &
    \int_{t_o}^{T} \int_{\reali^n}
    [
    (u - k) (\partial_t \phi
    +
    c \cdot \nabla \phi)
    +
    (b \,u - k \, \div c) \phi
    ] \sgn (u - k) \,
    \d{x} \, \d{t}
    \\
    \!\!\!\!\!\!
    & = &
    \int_{t_o}^{T} \int_{\reali^n}
    \Bigl[
    \left(
      \frac{u_o (y) \, \mathcal{B} (t,y)}{J (t,y)} - k
    \right)
    \left(
      \partial_t \phi \left(t,X (t; t_o,y)\right)
      +
      c \left(t,X (t; t_o,y) \right)
      \cdot
      \nabla \phi \left(t,X (t; t_o,y)\right)
    \right)
    \\
    \!\!\!\!\!\!
    & &
    +
    \left( b \left(t,X (t; t_o,y)\right) \,
      \frac{u_o (y) \, \mathcal{B} (t,y)}{J (t,y)}
      -
      k \, \div c \left(t,X (t; t_o,y)\right) \right)
    \phi \left(t,X (t; t_o,y)\right) \Bigr]
    \\
    \!\!\!\!\!\!
    & &
    \qquad \times \sgn \left( \frac{u_o (y) \, \mathcal{B} (t,y)}{J (t,y)} -
      k \right) \, J (t,y) \, \d{y} \, \d{t}
    \\
    \!\!\!\!\!\!
    & = &
    \int_{t_o}^{T} \int_{\reali^n} \Bigl[ u_o (y) \,
    \mathcal{B} (t,y) \,\frac{\d{}}{\d{t}} \phi \left(t,X (t; t_o,y)\right)
    -
    k \,
    J (t,y) \, \frac{\d{}}{\d{t}} \phi \left(t,X (t; t_o,y)\right)
    \\
    \!\!\!\!\!\!
    & &
    + \phi \left(t,X (t; t_o,y)\right) \, \frac{\d{}}{\d{t}} \left( u_o (y) \,
      \mathcal{B} (t,y) \right) - k \, \phi \left(t,X (t; t_o,y)\right) \,
    \frac{\d{}}{\d{t}} J(t,y) \Bigr]
    \\
    \!\!\!\!\!\!
    & &
    \qquad \times \sgn \left( u_o (y) \, \mathcal{B} (t,y) - k \, J (t,y)
    \right) \, \d{y} \, \d{t}
    \\
    \!\!\!\!\!\!
    & = & \int_{t_o}^{T} \int_{\reali^n} \frac{\d{}}{\d{t}} \left[
      \left( u_o (y) \, \mathcal{B} (t,y) - k \, J (t,y) \right) \, \phi
      \left(t,X (t; t_o,y)\right) \right]
    \\
    \!\!\!\!\!\!
    & &
    \qquad
    \times \sgn \left( u_o (y) \,\mathcal{B} (t,y) - k \, J (t,y)
    \right) \, \d{y} \, \d{t}
    \\
    \!\!\!\!\!\!
    & = &
    \int_{t_o}^{T} \int_{\reali^n} \frac{\d{}}{\d{t}}
    \left(
      \modulo{u_o (y) \, \mathcal{B} (t,y) - k \, J (t,y) } \,
      \phi \left(t,X (t; t_o,y)\right)
    \right)
    \, \d{y} \, \d{t}
    \\
    \!\!\!\!\!\!
    & = &
    0 \,.
  \end{eqnarray*}
  It is immediate to verify that for $u$ as in~\eqref{eq:solCara} the
  initial condition holds. By~\cite[Section~2]{Kruzkov}, $u$ is also a
  weak solution.

  \smallskip

\noindent\textbf{2.} Let $u$ be defined as in~\eqref{eq:solCara} and
$u_*$ be a weak solution to~\eqref{eq:CauchyH}. Then, by the step
above, the function $U = u - u_*$ is a weak solution
to~\eqref{eq:CauchyH} with zero initial datum.

Fix $\tau \in \left]t_o,T\right]$, choose any $\phi \in \C1 (I;\Cc1
(\reali^n;\reali))$ and let $\beta_\epsilon \in \C1 (I;\reali)$ such
that $\beta_\epsilon (t) = 1$ for $t \in [t_o + \epsilon,\tau -
\epsilon]$, $\beta_\epsilon' (t) \in [0,2/\epsilon]$ for $t \in
[t_o,t_o + \epsilon]$, $\beta_\epsilon' (t) \in [-2/\epsilon,0]$ for
$t \in [\tau - \epsilon,\tau]$ and $\beta_\epsilon (t) = 0$ for $t \in
[\tau,T]$. Using the definition of weak solution,
\begin{eqnarray*}
  0 & = &
  \int_I \int_{\reali^n}
  \left(
    U \, \partial_t (\phi\, \beta_\epsilon)
    +
    U c \cdot \nabla (\phi \, \beta_\epsilon)
    +
    b \, U \, \phi \, \beta_\epsilon
  \right)
  \d{x} \, \d{t}
  \\
  & = &
  \int_{t_o}^\tau \!\! \int_{\reali^n}
  U \left(
    \partial_t \phi
    +
    c \cdot \nabla \phi
    +
    b \, \phi \right) \beta_\epsilon \,
  \d{x} \, \d{t}
  +
  \int_{t_o}^{t_o + \epsilon} \!\!\!
  \beta_\epsilon'
  \int_{\reali^n} \!\!\!
  U \, \phi \, \d{x} \, \d{t}
  +
  \int_{\tau-\epsilon}^\tau \!\!\!
  \beta_\epsilon'
  \int_{\reali^n} \!\!\! U \, \phi \, \d{x} \, \d{t}
  \\
  & = &
  \int_{t_o}^\tau
  \beta_\epsilon
  \int_{\reali^n}
  U
  \left(
    \partial_t \phi + c \cdot \nabla \phi + b \, \phi
  \right) \, \d{x} \, \d{t}
  \\
  & & \quad
  +
  \int_0^1
  \epsilon \, \beta_\epsilon' (\epsilon t + t_o)
  \int_{\reali^n} U (\epsilon t + t_o,x) \, \phi (\epsilon t + t_o,x)\,
  \d{x} \, \d{t}
  \\
  & &
  \quad
  +
  \int_{\tau-\epsilon}^\tau \beta_\epsilon' (t)
  \int_{\reali^n} U (t,x)\, \phi (t,x)\, \d{x} \, \d{t} \,.
\end{eqnarray*}
As $\epsilon \to 0$, the first term converges to $\int_{t_o}^\tau
\int_{\reali^n} U \left( \partial_t \phi + c \cdot \nabla \phi + b \,
  \phi \right) \, \d{x} \, \d{t}$. By the Dominated Convergence
Theorem, the second term tends to $0$, since $U (t_o) = 0$. Concerning
the third term, note that
\begin{eqnarray*}
  & &
  \modulo{
    \int_{\tau-\epsilon}^\tau \beta_\epsilon' (t)
    \int_{\reali^n} U (t,x)\, \phi (t,x)\, \d{x} \, \d{t}
    +
    \int_{\reali^n} U (\tau,x) \, \phi (\tau,x) \, \d{x}
  }
  \\
  & = &
  \modulo{
    \int_{\tau-\epsilon}^\tau
    \beta_\epsilon' (t)
    \int_{\reali^n}
    \left(
      U (t,x)\, \phi (t,x) - U (\tau,x) \, \phi (\tau,x)
    \right)
    \, \d{x} \, \d{t} } \to 0
  \quad \mbox{ as } \quad \epsilon \to 0
\end{eqnarray*}
by the continuity of $U$ in time and the smoothness of $\phi$.

Choose now any $\eta \in \Cc1 (\reali^n;\reali)$ and define $\phi$ as
the backward solution to
\begin{equation*}
  \left\{
    \begin{array}{l}
      \partial_t \phi + c \cdot \nabla \phi =- b\, \phi \\
      \phi (\tau) = \eta \,.
    \end{array}
  \right.
\end{equation*}
Then,
\begin{displaymath}
  0
  =
  \int_{t_o}^\tau \int_{\reali^n}
  U
  \left( \partial_t \phi + c \cdot \nabla \phi + b \, \phi \right)
  \, \d{x} \, \d{t}
  -
  \int_{\reali^n} U (\tau,x) \, \phi (\tau, x) \, \d{x}
  =
  - \int_{\reali^n} U (\tau,x) \, \eta (x) \, \d{x}
\end{displaymath}
proving that $U (\tau)$ vanishes identically. By the arbitrariness of
$\tau$, $u=u_*$.
\end{proofof}

\begin{proofof}{Proposition~\ref{prop:hyper}}
  Define $g (t,x) = b (t,x) - \div c (t,x)$. The equation
  in~\eqref{eq:CauchyH} fits into the general form
  \begin{displaymath}
    \partial_t u + \div f (t,x,u) = F (t,x,u)
    \quad
    \mbox{ with } \quad
    \begin{array}{rcl}
      f (t,x, u) & = & c (t,x) \, u
      \\
      F (t,x,u) & = & b (t, x) \, u \,.
    \end{array}
  \end{displaymath}
  Points~\textbf{1.} and~\textbf{2.} follow, for instance,
  from~\cite[Theorem~1 and Theorem~2]{Kruzkov}.

  \noindent\textbf{3.}\quad
  Start from~\cite[Formula~(3.1)]{Kruzkov} and use the fact that $0$
  solves the equation in~\eqref{eq:CauchyH}, getting
  \begin{displaymath}
    \norma{\mathcal{H}_{t_o,t} u_o}_{\L1 (\reali^n; \reali)}
    \leq
    \norma{u_o}_{\L1 (\reali^n; \reali)} \,
    \exp\left(\norma{b}_{\L\infty ([t_o,t]\times \reali^n; \reali)} (t-t_o)\right)
    \,.
  \end{displaymath}

  \noindent\textbf{4.}\quad From~\eqref{eq:solCara} it easily follows
  that
  \begin{displaymath}
    \norma{\mathcal{H}_{t_o,t} u_o}_{\L\infty (\reali^n; \reali)}
    \leq
    \norma{u_o}_{\L\infty (\reali^n; \reali)} \,
    \exp\left(\norma{g}_{\L\infty([t_o,t]\times \reali^n; \reali)} (t-t_o)\right)
    \,.
  \end{displaymath}

  \noindent\textbf{5.}\quad We refer
  to~\cite[Theorem~2.6]{ColomboMercierRosini}, see
  also~\cite{ColomboMercierRosiniCR}, as refined
  in~\cite[Proposition~2.10]{MercierDaSola}
  and~\cite[Proposition~2.9]{MercierV2}. Indeed, compute preliminarily
  \begin{equation}
    \label{eq:table}
    \begin{array}{c}
      \begin{array}{rcl@{\qquad\qquad\qquad}rcl}
        \partial_u f (t,x,u) & = & c (t,x)
        &
        \div f (t,x,u) & = & \left(\div c (t,x)\right) \; u
        \\
        \partial_u \nabla f (t,x,u) & = & \nabla c (t,x)
        &
        \nabla^2 f (t,x,u) & = & \left(\nabla^2 c (t,x)\right) \; u
        \\
        \partial_u F (t,x,u) & = & b (t,x)
        &
        \nabla F (t,x,u) & = & \left(\nabla b (t,x)\right) \; u
      \end{array}
      \\
      \begin{array}{rclcl}
        F (t,x,u) - \div f (t,x,u)
        & = &
        \left(b (t,x) - \div c (t,x)\right) u
        & = &
        g (t,x) \, u
        \\
        \partial_u \left(F (t,x,u) - \div f (t,x,u)\right)
        & = &
        b (t,x) - \div c (t,x)
        & = &
        g (t,x)
        \\
        \nabla\left(F (t,x,u) - \div f (t,x,u)\right)
        & = &
        \left(\nabla b (t,x) - \nabla \left(\div c (t,x)\right)\right) u
        & = &
        \nabla g (t,x) \, u
      \end{array}
    \end{array}
  \end{equation}
  and introduce the quantities
  \begin{displaymath}
    \begin{array}{@{}r@{\,}c@{\,}l@{}}
      \kappa_o^*
      & = &
      (2n+1) \norma{\nabla c_1}_{\L\infty ([t_o,t] \times \reali^n; \reali^{n \times n})}
      +
      \norma{b_1}_{\L\infty ([t_o,t] \times \reali^n; \reali)}
      \\
      \kappa^*
      & = &
      \norma{b_1}_{\L\infty([t_o,t] \times \reali^n; \reali)}
      \!+\!
      \norma{\div (c_1-c_2)}_{\L\infty([t_o,t] \times \reali^n; \reali)}
      \\
      \kappa^*_1
      & = &
      \norma{b_1}_{\L\infty([t_o,t] \times \reali^n; \reali)}
      \!+\!
      \norma{\div (c_1-c_2)}_{\L\infty([t_o,t] \times \reali^n; \reali)}
      \!+\!
      (2n+1) \norma{\nabla c_1}_{\L\infty ([t_o,t] \times \reali^n; \reali^{n \times n})}.
    \end{array}
  \end{displaymath}
  For $i=1,2$, let $b_i$ satisfy~\textbf{(b)} with $b_1-b_2 \in \L1 (I
  \times \reali^n; \reali)$ and $c_i$ satisfy~\textbf{(c)} with $\div
  (c_1 -c_2) \in \L1 (I \times \reali^n; \reali)$. Then, it is
  immediate to check that the requirements
  in~\cite[Section~2]{MercierDaSola} and
  in~\cite[Section~2]{MercierV2} hold. Indeed, with obvious notation,
  \begin{eqnarray*}
    \left.
      \begin{array}{r}
        c_i \in \L\infty (I \times \reali^n; \reali^n); \;
        c_i \in \C2 (I \times \reali^n; \reali^n); \;
        b_i \in \C1(I \times \reali^n; \reali)
        \\
        (b_i - \div c_i) \in
        \L\infty (I \times \reali^n; \reali)
      \end{array}
    \right\}
    & \Rightarrow&
    \mbox{\textbf{(H1*)}}
    \\
    \left.
      \begin{array}{r}
        \nabla c_1 \in
        \L\infty (I \times \reali^n; \reali^{n\times n}); \;
        b_1 \in \L\infty (I \times \reali^n; \reali); \;
        \\
        \nabla b_1 - \nabla (\div c_1) \in \L1 (I \times \reali^n; \reali^n)
      \end{array}
    \right\}
    & \Rightarrow&
    \mbox{\textbf{(H2*)}}
    \\
    \left.
      \begin{array}{r}
        c_1-c_2 \in \L\infty (I\times \reali^n; \reali^n)
        \\
        (b_1-b_2) \in \L\infty (I \times \reali; \reali); \;
        \left((b_1-b_2) - \div (c_1-c_2)\right) \in
        \L1 (I \times \reali^n; \reali)
      \end{array}
    \right\}
    & \Rightarrow&
    \mbox{\textbf{(H3*)}}
  \end{eqnarray*}
  Note also that
  \begin{displaymath}
    \frac{\kappa_o^* e^{\kappa_o^*t} - \kappa^* e^{\kappa*t}}{\kappa_o^* - \kappa^*}
    =
    e^{\kappa_o^* t}
    +
    \frac{\kappa^* \left(e^{\kappa_o^*t} - e^{\kappa^* t}\right)}{\kappa_o^* - \kappa^*}
    \leq
    (1+\kappa^*\,t) e^{\kappa_1^* t}
    \leq
    e^{\kappa^* t} e^{\kappa_1^* t}.
  \end{displaymath}
  Applying~\cite[Proposition~2.10]{MercierDaSola}
  and~\cite[Proposition~2.9]{MercierV2} we now obtain
  \begin{eqnarray*}
    & &
    \norma{
      \mathcal{H}^1_{t_o,t} u_o - \mathcal{H}^2_{t_o,t} u_o}_{\L1 (\reali^n; \reali)}
    \\
    & \leq &
    e^{ (\kappa^*+\kappa_1^*) \, (t-t_o)} \,
    \norma{c_1-c_2}_{\L1 ([t_o,t]; \L\infty (\reali^n; \reali^n))} \,
    \\
    & &
    \times
    \Bigl[
    \tv (u_0)
    +
    I_n
    \int_{t_o}^t
    e^{-\kappa_o^* (\tau - t_o)}
    \norma{\mathcal{H}^1_{t_o,\tau} u_o}_{\L\infty (\reali^n;\reali)}
    \norma{\nabla g_1 (\tau)}_{\L1 (\reali^n; \reali^n)}
    \d\tau
    \Bigr]
    \\
    & &
    +
    e^{\kappa^* (t - t_o)}
    \int_{t_o}^t
    \norma{g_1 (\tau) - g_2 (\tau)}_{\L{1} (\reali^n; \reali)}
    \max_{i=1,2}
    \norma{\mathcal{H}^i_{t_o,\tau} u_o}_{\L\infty (\reali^n;\reali)}
    \d\tau
  \end{eqnarray*}
  where $\norma{\mathcal{H}^i_{t_o,t} (u_o)}_{\L\infty
    (\reali^n;\reali)}$ is estimated in~4., $g_i = b_i - \div c_i$,
  $I_n$ is as in~\eqref{eq:3H}.

  \noindent\textbf{6.}\quad Directly follows from~\eqref{eq:solCara}.

  \noindent\textbf{7.}\quad By assumptions, $b$ satisfies~\textbf{(b)} and $c$ satisfies~\textbf{(c)}, hence by~\eqref{eq:table}, both~\textbf{(H1*)} and~\textbf{(H2*)} hold. From~\cite[Theorem~2.5]{ColomboMercierRosini}
  or~\cite[Theorem~2.2]{MercierV2}, it directly follows that
  \begin{equation}
    \label{eq:TV}
    \!\!\!
    \tv\left(\mathcal{H}_{t_o,t} u_o \right)
    \leq
    \tv (u_o) e^{\kappa_o^* (t-t_o)}
    +
    I_n \!
    \int_{t_o}^t \!
    e^{\kappa_o^* (t-\tau)}
    \norma{\nabla g (\tau)}_{\L1 (\reali^n;\reali^n)}
    \norma{\mathcal{H}_{t_o,\tau} u_o}_{\L\infty (\reali^n; \reali)}
    \d\tau
  \end{equation}
  where $\norma{\mathcal{H}_{t_o,t} (u_o)}_{\L\infty
    (\reali^n;\reali)}$ is estimated in~4., $I_n$ is as
  in~\eqref{eq:3H} and
  \begin{displaymath}
    \kappa_o^*
    =
    (2n+1) \norma{\nabla c}_{\L\infty ([t_o,t] \times \reali^n; \reali^{n \times n})}
    +
    \norma{b}_{\L\infty ([t_o,t] \times \reali^n; \reali)} \,.
  \end{displaymath}

  \noindent\textbf{8.}\quad Assume that $t_1 < t_2$. From
  Definition~\ref{def:Hyp}, we have that
  \begin{equation}
    \label{eq:ti}
    \int_{t_i}^T \int_{\reali^n}
    \left(
      u \, \partial_t \phi + u \, c \cdot \nabla \phi + b \, u \, \phi
    \right)
    \d{x}\d{t}
    +
    \int_{\reali^n} \phi (t_i, x) \, u (t_i,x) \d{x}
    =
    0 \,.
  \end{equation}
  Following the proof of~\cite[Theorem~4.3.1]{DafermosBook}, let $\phi
  (t,x) = \chi (t) \, \psi (x)$ with $\chi \in \Cc\infty (\pint{I};
  \reali)$, $\chi (t) = 1$ for all $t \in [t_1, t_2]$ and $\psi \in
  \Cc\infty (\reali^n; \reali)$ with $\modulo{\psi (x)} \leq 1$ for
  all $x \in \reali^n$. Subtract~\eqref{eq:ti} for $i=1$
  from~\eqref{eq:ti} for $i=2$,
  use~\cite[Proposition~3.2]{AmbrosioFuscoPallara} and the estimates
  at points~3., 4.~and~7.~to obtain:
  \begin{align*}
    & \norma{u (t_2)-u (t_1)}_{\L1 (\reali^n; \reali)}
    \\
    = {} & \sup_{\psi \in \Cc\infty, \,\modulo{\psi (x)} \leq 1}
    \int_{\reali^n} \psi (x) \left(u (t_2,x) - u (t_1,x)\right) \d{x}
    \\
    = {} & \sup_{\psi \in \Cc\infty, \,\modulo{\psi (x)} \leq 1}
    \int_{t_1}^{t_2} \!\!\! \int_{\reali^n} \!\!  \left( u (t,x)\, c
      (t,x) \cdot \nabla \psi (x) + b (t,x) \, u (t,x) \, \psi
      (x)\right) \d{x} \d{t}
    \\
    = {} & \sup_{\psi \in \Cc\infty, \,\modulo{\psi (x)} \leq 1}
    \int_{t_1}^{t_2} \!\!\! \int_{\reali^n} \!\!  \left[ -\nabla \!\!
      \left(u (t,x) c (t,x)\right) \psi (x) + b (t,x) u (t,x) \psi (x)
    \right] \d{x} \d{t}
    \\
    \leq {} & \int_{t_1}^{t_2} \!\!\! \int_{\reali^n} \!\!  \left[
      \norma{\nabla \!\! \left(u (t,x) c (t,x)\right)} + \modulo{b
        (t,x) u (t,x)} \right] \d{t} \d{x}
    \\
    \leq {} & \int_{t_1}^{t_2} \!\!  \left[ \int_{\reali^n} \norma{u
        (t,x) \nabla c (t,x)} \d{x} + \int_{\reali^n} \norma{c (t,x)}
      \d{\left(\nabla u (t)\right) (x)} + \int_{\reali^n} \modulo{b
        (t,x) u (t,x)} \d{x} \right] \d{t}
    \\
    \leq {} & \int_{t_1}^{t_2} \norma{u (t)}_{\L1 (\reali^n; \reali)}
    \left[ \norma{b (t)}_{\L\infty (\reali^n; \reali)} + \norma{\nabla
        c (t)}_{\L\infty (\reali^n; \reali^{n\times n})} \right] \d{t}
    \\
    & + \int_{t_1}^{t_2} \norma{c (t)}_{\L\infty (\reali^n; \reali^n)}
    \tv\left(u (t) \right) \d{t}
    \\
    \leq {} & \norma{u_o}_{\L1 (\reali^n; \reali)}
    e^{\norma{b}_{\L\infty ([t_o,t_2]\times \reali^n; \reali)}
      (t_2-t_o)} \left[ \norma{b}_{\L\infty ([t_1,t_2]\times\reali^n;
        \reali)} + \norma{\nabla c}_{\L\infty
        ([t_1,t_2]\times\reali^n; \reali^{n\times n})} \right]
    \modulo{t_2-t_1}
    \\
    & + \norma{c}_{\L\infty ([t_1,t_2]\times \reali^n; \reali^n)}
    \int_{t_1}^{t_2} \tv (u_o) e^{\kappa_o^* (t-t_o)} \d{t}
    \\
    & + I_n \norma{c}_{\L\infty ([t_1,t_2]\times \reali^n; \reali^n)}
    \int_{t_1}^{t_2} \int_{t_o}^t\!  e^{\kappa_o^* (t-\tau)} \norma{
      \nabla g (\tau)}_{\L1 (\reali^n;\reali^n)} \norma{u
      (\tau)}_{\L\infty (\reali^n; \reali)} \!  \d\tau \d{t}
    \\
    \leq {} & \modulo{t_2-t_1} \Bigg[ \norma{u_o}_{\L1 (\reali^n;
      \reali)} e^{\norma{b}_{\L\infty ([t_o,t_2]\times \reali^n;
        \reali)} (t_2-t_o)} \left[ \norma{b}_{\L\infty
        ([t_1,t_2]\times\reali^n; \reali)} + \norma{\nabla
        c}_{\L\infty ([t_1,t_2]\times\reali^n; \reali^{n\times n})}
    \right]
    \\
    & + \norma{c}_{\L\infty ([t_1,t_2]\times \reali^n; \reali^n)} \tv
    (u_o) e^{\kappa_o^* (t_2-t_o)}
    \\
    & +\!  I_n \norma{u_o}_{\L\infty (\reali^n; \reali)}
    \norma{c}_{\L\infty ([t_1,t_2]\times \reali^n; \reali^n)}\!\!
    \int_{t_o}^{t_2}\!\!\!  \norma{\nabla g (\tau)}_{\L1 (\reali^n;
      \reali^n)} e^{\kappa_o^*
      (t_2-\tau)+\norma{g}_{\L\infty([t_o,\tau]\times\reali^n;
        \reali)}(\tau-t_o)}\!  \d\tau\!  \Bigg]
  \end{align*}
  proving point~8.

  \noindent\textbf{9.}\quad Directly follows from~\eqref{eq:solCara},
  since, by~\textbf{(c)}, the speed of characteristics is bounded.
\end{proofof}

\subsection{Proof of Theorem~\ref{thm:main}}
\label{subs:Main}

\begin{proofof}{Theorem~\ref{thm:main}}
  We are going to construct a solution to~\eqref{eq:Model} as limit of
  a Cauchy sequence of approximate solutions in the complete metric
  space
  \begin{displaymath}
    X = \L1([0,T]; \mathcal{X}^+),
  \end{displaymath}
  for a suitable positive $T$ to be chosen below, equipped with the
  distance
  \begin{eqnarray*}
    d\left((u_1,w_1), (u_2, w_2)\right)
    & = &
    \norma{u_2-u_1}_{\L1 ([0,T];\L1 (\reali^n;\reali))}
    +
    \norma{w_2-w_1}_{\L1 ([0,T];\L1 (\reali^n;\reali))}
    \\
    & = &
    \int_0^T \int_{\reali^n}
    \left(
      \modulo{u_2 (t,x) - u_1 (t,x)}
      +
      \modulo{w_2 (t,x) - w_1 (t,x)}
    \right)
    \d{x} \d{t} \,.
  \end{eqnarray*}
  For $r>0$, we introduce the domain
  \begin{displaymath}
    \mathcal{X}_r
    =
    \left\{
      (u,w) \in \mathcal{X}^+ \colon
      \begin{array}{rcl@{,\qquad}rcl}
        \norma{u}_{\L\infty (\reali^n;\reali)} & \leq & r
        &
        \tv (u) & \leq & r
        \\
        \norma{w}_{\L\infty (\reali^n;\reali)} & \leq & r
        &
        \norma{w}_{\L1 (\reali^n;\reali)} & \leq & r
      \end{array}
    \right\}.
  \end{displaymath}
  Choose an initial datum
  \begin{equation}
    \label{eq:woRegular}
    (u_o, w_o) \in \mathcal{X}_r
    \quad \mbox{ with moreover } \quad
    w_o \in (\C1 \cap\W11) (\reali^n; \reali^+) \,,
  \end{equation}
  set, for $t \in [0, T]$, $\left(u_0 (t),w_0 (t)\right) =
  (u_o,w_o)$. For $i \in \naturali$, define recursively for $(t,x) \in
  [0,T] \times \reali^n$,
  \begin{displaymath}
    a_{i+1} (t,x) = \gamma - \delta \, u_i (t,x)\,,
    \quad
    b_{i+1} (t,x) = \alpha \, w_i (t,x) - \beta
    \quad \mbox{ and } \quad
    c_{i+1} (t,x) = v \left(w_i (t)\right) (x)
  \end{displaymath}
  and let $(u_{i+1}, w_{i+1})$ be such that
  \begin{displaymath}
    \left\{
      \begin{array}{l}
        \partial_t u_{i+1} + \div \left(c_{i+1} u_{i+1}\right) = b_{i+1} \, u_{i+1}
        \\
        \partial_t w_{i+1} - \mu \, \Delta w_{i+1} = a_{i+1} \, w_{i+1}
        \\
        u_{i+1} (0) = u_o
        \\
        w_{i+1} (0) = w_o \,.
      \end{array}
    \right.
  \end{displaymath}

  \bigskip

  \noindent\textbf{Claim~0:} For all $i \in \naturali$,
  \begin{description}
    \setlength{\itemsep}{0pt} \setlength{\parskip}{0pt}
  \item[C0.1] $(u_i, w_i)$ is well defined and in $X$;
  \item[C0.2] $w_i \in \C1 ([0,T] \times \reali^n ; \reali)$ and
    $\nabla w_i \in \L1 ([0,T] \times \reali^n; \reali^n)$;
  \item[C0.3] $a_{i+1} \in \L\infty ([0,T] \times \reali^n; \reali)$;
  \item[C0.4] $b_{i+1}$ satisfies~\textbf{(b)} with $I=[0,T]$;
  \item[C0.5] $c_{i+1}$ satisfies~\textbf{(c)} with $I=[0,T]$.
  \end{description}

  \noindent\textbf{Proof of Claim~0.}  We prove it by induction.

  \noindent \textbf{Case $i=0$:} is immediate by~\eqref{eq:woRegular}
  and by the above definition of $a_1, b_1, c_1$.

  \noindent \textbf{From $i-1$ to $i$:} Assume now that~C0.1, $\ldots$
  ,~C0.5 are all satisfied up to the $i$-th iteration. Then,
  Proposition~\ref{prop:para} and Proposition~\ref{prop:hyper} can now
  be applied, proving~C0.1. Moreover, by~7.~and
  9.~in~Proposition~\ref{prop:para}, also~C0.2 holds. Furthermore, the
  estimate at~3.~and~4.~in Proposition~\ref{prop:hyper} ensure
  that~C0.3 holds. Moreover, C0.1 and C0.2 directly imply~C0.4 and,
  together with~\vv, also~C0.5, completing the proof of the present
  claim.

  \medskip

  From the above, thanks to~{\vv}, it clearly follows that:
  \begin{equation}
    \label{eq:1}
    a_{i+1} - a_1 \in \L1 (I\times \reali^n; \reali)
    \,,\;
    b_{i+1} - b_i \in \L1 (I \times \reali^n; \reali)
    \,,\;
    \div (c_{i+1} - c_i) \in \L1 (I\times \reali^n; \reali) \,.
  \end{equation}

  In the next two claims we particularize the $\L1$ and $\L\infty$
  estimates in~3.~and~4.~of Propositions~\ref{prop:para}
  and~\ref{prop:hyper}, thanks to the explicit expressions of $a$, $b$
  and $c$.

  \medskip

  \noindent\textbf{Claim~1:} For all $i \in \naturali$, if $w_i$ is
  defined up to time $\hat T$, then for all $t \in [0, \hat T]$,
  \begin{displaymath}
    \norma{ w_i (t) }_{\L1 (\reali^n; \reali)}
    \leq
    \norma{ w_o }_{\L1 (\reali^n; \reali)} \, e^{\gamma t}
    \quad \mbox{and} \quad
    \norma{ w_i (t) }_{\L\infty (\reali^n; \reali)}
    \leq
    \norma{ w_o }_{\L\infty (\reali^n; \reali)} \, e^{\gamma t}.
  \end{displaymath}

  \noindent\textbf{Proof of Claim~1.} Assume $i>0$, the case $i=0$
  being obvious. By~\eqref{eq:solNucleo},
  \begin{eqnarray*}
    w_i (t,x)
    & = &
    \int_{\reali^n} \!\!\! H_\mu (t, x-\xi) \, w_o (\xi) \d\xi
    +
    \int_0^t \!\! \int_{\reali^n} \!\!\! H_\mu (t - \tau, x-\xi)
    \left(\gamma - \delta \, u_{i-1} (\tau, \xi)\right)
    w_i (\tau, \xi) \d\xi \d\tau
    \\
    & \leq &
    \int_{\reali^n} H_\mu (t, x-\xi) \, w_o (\xi) \, \d\xi +
    \int_0^t \int_{\reali^n} \gamma \, H_\mu (t - \tau, x-\xi) \, w_i
    (\tau, \xi) \d\xi \d\tau.
  \end{eqnarray*}
  By Gronwall Lemma and~\eqref{eq:calore1}:
  \begin{eqnarray*}
    w_i (t,x)
    & \leq &
    \left(
      \int_{\reali^n}
      H_\mu (t, x-\xi) \, w_o (\xi) \,
      \d\xi
    \right)
    \exp\left(
      \int_0^t \int_{\reali^n} \gamma \,
      H_\mu (t - \tau, x-\xi) \, \d\xi \, \d\tau
    \right)
    \\
    & = &
    e^{\gamma t} \, \int_{\reali^n} H_\mu (t, x-\xi) \, w_o (\xi)
    \, \d\xi.
  \end{eqnarray*}
  The proof of the claim follows.

  \bigskip

  \noindent\textbf{Claim~2:} For all $i \in \naturali$, if $u_i$ is
  defined up to time $\hat T$, then for all $t \in [0,\hat T]$,
  \begin{eqnarray*}
    \norma{u_i (t)}_{\L1 (\reali^n; \reali)}
    & \leq &
    \norma{u_o}_{\L1 (\reali^n; \reali)} \,
    \exp \left(
      \alpha \, \frac{e^{\gamma t}-1}{\gamma} \,
      \norma{w_o}_{\L\infty (\reali^n;\reali)}
    \right)
    \\
    & \leq &
    \norma{u_o}_{\L1 (\reali^n; \reali)} \,
    \exp \left(
      \alpha \, t \, e^{\gamma t} \,
      \norma{w_o}_{\L\infty (\reali^n;\reali)}
    \right)
    \\
    \norma{u_i (t)}_{\L\infty (\reali^n; \reali)}
    & \leq &
    \norma{u_o}_{\L\infty (\reali^n; \reali)} \,
    \exp \left(
      (\alpha+K) \, \frac{e^{\gamma t}-1}{\gamma} \,
      \norma{w_o}_{\L\infty (\reali^n;\reali)}
    \right)
\\
    & \leq &
    \norma{u_o}_{\L\infty (\reali^n; \reali)} \,
    \exp \left(
      (\alpha+K) \, t \, e^{\gamma t} \,
      \norma{w_o}_{\L\infty (\reali^n;\reali)}
    \right).
  \end{eqnarray*}

  \noindent\textbf{Proof of Claim~2.} Assume $i>0$, the case $i=0$
  being obvious. By~{\vv} and Claim~0, we can apply
  Lemma~\ref{lem:hypSol} and by~\eqref{eq:solCara} we obtain
  \begin{equation}
    \label{eq:ineq}
    u_i (t,x)
    \leq
    u_o\left(X (0;t,x)\right) \,
    \exp
    \int_0^t \left(
      \alpha \, w_{i-1}\left(\tau, X (\tau; t,x)\right)
      -
      \div c_i\left(\tau, X (\tau; t,x)
      \right) \d\tau
    \right) \,.
  \end{equation}
  To obtain the $\L1$ estimate, we adopt the notation
  in~\eqref{eq:CambioVar} with $t_o=0$, $b = \alpha w_{i-1}$ and $c =
  c_i$, so that, using Claim~1 above, we have
  \begin{eqnarray*}
    \norma{u_i (t)}_{\L1 (\reali^n; \reali)}
    & \leq &
    \int_0^t
    \frac{1}{J (\tau,y)} \, u_o (y) \, \mathcal{B} (\tau,y) \, J (\tau,y) \,
    \d\tau
    \\
    & \leq &
    \norma{u_o}_{\L1 (\reali^n; \reali)} \,
    \exp\left(
      \int_0^t \alpha \,
      \norma{w_o}_{\L\infty (\reali^n;\reali)} \, e^{\gamma\tau} \, \d\tau
    \right)
    \\
    & = &
    \norma{u_o}_{\L1 (\reali^n; \reali)} \,
    \exp \left(
      \alpha \, \frac{e^{\gamma t}-1}{\gamma} \,
      \norma{w_o}_{\L\infty (\reali^n;\reali)}
    \right) .
  \end{eqnarray*}
  The $\L\infty$ estimate is obtained from~\eqref{eq:ineq} using~{\vv}
  and Claim~1:
  \begin{eqnarray*}
    \norma{u_i (t)}_{\L\infty (\reali^n;\reali)}
    & \! \leq &
    \norma{u_o}_{\L\infty (\reali^n;\reali)} \,
    \exp
    \int_0^t
    \left(
      \alpha \, \norma{w_{i-1} (\tau)}_{\L\infty (\reali^n;\reali)}
      +
      K \, \norma{w_{i-1} (\tau)}_{\L\infty (\reali^n; \reali)}
    \right)
    \d\tau
    \\
    & \! \leq &
    \norma{u_o}_{\L\infty (\reali^n;\reali)} \,
    \exp
    \int_0^t
    (\alpha+K)  \, \norma{w_o}_{\L\infty (\reali^n;\reali)} \, e^{\gamma \tau} \,
    \d\tau
    \\
    & \!= &
    \norma{u_o}_{\L\infty (\reali^n; \reali)} \,
    \exp \left(
      (\alpha+K) \, \frac{e^{\gamma t}-1}{\gamma} \,
      \norma{w_o}_{\L\infty (\reali^n;\reali)}
    \right),
  \end{eqnarray*}
  completing the proof of Claim~2.

  \bigskip

  We now prove that there exist positive $T$ and $\mathcal{K} (T,r)$
  such that for all $i \in \naturali$,
  \begin{equation}
    \label{eq:distanza}
    d\left((u_{i+1}, w_{i+1}), (u_i, w_i)\right)
    \leq
    T \, \mathcal{K} (T,r) \,
    d\left((u_i, w_i), (u_{i-1}, w_{i-1})\right).
  \end{equation}
  By~\eqref{eq:1}, recall the proof of~5.~in
  Proposition~\ref{prop:para} and apply the $\L\infty$ estimate in
  Claim~1:
  \begin{eqnarray*}
    \norma{w_{i+1} (t) - w_i (t)}_{\L1 (\reali^n; \reali)}
    & \leq &
    \int_0^t \norma{a_{i+1} (\tau) - a_i (\tau)}_{\L1 (\reali^n;
      \reali)} \norma{w_{i+1} (\tau)}_{\L\infty (\reali^n; \reali)} \,
    \d\tau
    \\
    & & + \int_0^t \norma{a_i (\tau)}_{\L\infty (\reali^n; \reali)}
    \norma{w_{i+1} (\tau) - w_i (\tau)}_{\L1 (\reali^n; \reali)} \,
    \d\tau
    \\
    & \leq & \norma{w_o}_{\L\infty (\reali^n; \reali)} \int_0^t
    e^{\gamma \, \tau} \, \norma{a_{i+1} (\tau) - a_i (\tau)}_{\L1
      (\reali^n; \reali)} \, \d\tau
    \\
    & &+ \int_0^t \norma{a_i (\tau)}_{\L\infty (\reali^n; \reali)}
    \norma{w_{i+1} (\tau) - w_i (\tau)}_{\L1 (\reali^n; \reali)} \,
    \d\tau.
  \end{eqnarray*}
  Apply Gronwall Lemma:
  \begin{eqnarray*}
    & &
    \norma{w_{i+1} (t) - w_i (t)}_{\L1 (\reali^n; \reali)}
    \\
    & \leq &
    \norma{w_o}_{\L\infty (\reali^n; \reali)} \,
    e^{\gamma t} \,
    \norma{a_{i+1} - a_i }_{\L1([0,t]\times \reali^n; \reali)} \,
    \exp \left(
      \norma{a_i}_{\L1 ([0,t]; \L\infty (\reali^n; \reali))} \right)
    \\
    & \leq &
    \norma{w_o}_{\L\infty (\reali^n; \reali)}\,
    e^{\gamma t} \,
    \norma{a_{i+1} - a_i }_{\L1([0,t] \times\reali^n; \reali)} \,
    \exp \left(
      t \, \norma{a_i}_{\L\infty([0,t] \times \reali^n; \reali)}
    \right).
  \end{eqnarray*}
  Since $ a_i (\tau,x) = \gamma - \delta \, u_{i-1} (\tau,x)$, by
  Claim~2
  \begin{eqnarray}
    \nonumber
    \norma{a_i}_{\L\infty ([0,t] \times \reali^n; \reali)}
    & \leq &
    \!\gamma
    +
    \delta \,
    \norma{u_{i-1} }_{\L\infty([0,t] \times \reali^n; \reali)}
    \\
    & \leq &
    \nonumber
    \!\gamma
    +
    \delta \, \norma{u_o}_{\L\infty(\reali^n;\reali)}
    \exp \!\left(\!
      (\alpha + K) \,
      \frac{e^{\gamma t}-1}\gamma \,
      \norma{w_o}_{\L\infty (\reali^n; \reali)} \!
    \right)\!.
    \\
    & \leq &
    \label{eq:an}
    \!\gamma
    +
    \delta \, \norma{u_o}_{\L\infty(\reali^n;\reali)}
    \exp \!\left(\!
      (\alpha + K) \,
      \norma{w_o}_{\L\infty (\reali^n; \reali)} \,
      t \, e^{\gamma t} \,
    \right)\!.
  \end{eqnarray}
  Hence,
  \begin{eqnarray}
    \nonumber
    & &
    \norma{w_{i+1} (t) - w_i (t)}_{\L1 (\reali^n; \reali)}
    \\
    \nonumber
    & \leq &
    \norma{w_o}_{\L\infty (\reali^n; \reali)} \,
    \delta \,
    \norma{u_i - u_{i-1} }_{\L1([0,t];\L1 (\reali^n; \reali))}
    \\
    \nonumber
    & &
    \times
    \exp \left(
      2 \, \gamma \, t
      +
      \delta \, t\,
      \norma{u_o}_{\L\infty(\reali^n; \reali)}
      \exp \left(
        (\alpha + K) \,
        \frac{e^{\gamma t}-1}\gamma \,
        \norma{w_o}_{\L\infty (\reali^n; \reali)}
      \right)
    \right) ;
    \\
    & & \nonumber
    \\
    \nonumber
    & &
    \norma{w_{i+1} - w_i }_{\L1 ([0,T]; \L1 (\reali^n; \reali))}
    \\
    \nonumber
    & \leq &
    \delta \, T \, \norma{w_o}_{\L\infty (\reali^n; \reali)}
    \, \norma{u_i - u_{i-1} }_{\L1([0,T];\L1 (\reali^n; \reali))}
    \\
    \nonumber
    & &
    \times
    \exp \left(
      2 \, \gamma \, T
      +
      \delta \, T\,
      \norma{u_o}_{\L\infty(\reali^n; \reali)}
      \exp \left(
        (\alpha + K) \,
        \frac{e^{\gamma T}-1}\gamma \,
        \norma{w_o}_{\L\infty (\reali^n; \reali)}
      \right)
    \right)
    \\
    \label{eq:diffw}
    & \leq &
    T \,
    \mathcal{K}_w (T,r) \, \norma{u_i - u_{i-1}}_{\L1([0,T];\L1 (\reali^n; \reali))},
  \end{eqnarray}
  since $(u_o,w_o) \in \mathcal{X}_r$, where
  \begin{equation}
    \label{eq:KTW}
    \mathcal{K}_w (T,r)
    =
    \delta \,
    r \,
    \exp \!
    \left[
      2 \, \gamma \, T
      +
      \delta \, T \,
      r \,
      \exp \left(
        (\alpha + K) \,
        \frac{e^{\gamma T}-1}\gamma \,
        r
      \right)
    \right] .
  \end{equation}


  We now pass to estimate $\norma{u_{i+1} - u_i}_{\L1 ([0,T];\L1
    (\reali^n;\reali))}$. To this aim, by~\eqref{eq:1}, we start from
  5.~in Proposition~\ref{prop:hyper} and use Claim~2 above:
  \begin{eqnarray*}
    & &
    \norma{
      u_{i+1} (t) - u_i (t)}_{\L1 (\reali^n; \reali)}
    \\
    & \leq &
    e^{ \kappa^* \, t} \,
    e^{\kappa_1^* \, t} \,
    \norma{c_{i+1}-c_i}_{\L1 ([0,t]; \L\infty (\reali^n; \reali^n))} \,
    \\
    & &
    \times \Bigl[
    \tv (u_0)
    +
    I_n \!
    \int_{0}^t \!
    e^{-\kappa_o^* \tau}
    \norma{u_{i+1} (\tau)}_{\L\infty (\reali^n;\reali)}
    \,
    \norma{
      \nabla b_{i+1} (\tau)
      -
      \nabla\left(\div c_{i+1} (\tau)\right)
    }_{\L1 (\reali^n; \reali^n)}
    \d\tau \Bigr]
    \\
    & &
    +
    \, e^{\kappa^* t} \,
    \int_{0}^t
    \norma{
      b_{i+1} (\tau) - b_i (\tau)
      -
      \div\left(c_{i+1} (\tau) - c_i (\tau)\right)
    }_{\L{1} (\reali^n; \reali)}
    \\
    & &
    \qquad\qquad\qquad
    \times \max
    \left\{
      \norma{u_{i+1} (\tau)}_{\L\infty (\reali^n;\reali)}
      ,
      \norma{u_i (\tau)}_{\L\infty (\reali^n;\reali)}
    \right\}
    \d\tau
    \\
    & \leq &
    e^{ \kappa^* \, t} \,
    e^{\kappa_1^* \, t} \,
    \norma{c_{i+1}-c_i}_{\L1 ([0,t]; \L\infty (\reali^n; \reali^n))} \,
    \Bigl[
    \tv (u_o)
    +
    I_n \, \norma{u_o}_{\L\infty (\reali^n;\reali)}
    \\
    & &
    \qquad
    \times
    \int_0^t
    e^{(\alpha+K) \norma{w_o}_{\L\infty (\reali^n; \reali)}\tau e^{\gamma \tau}}
    \norma{
      \nabla b_{i+1} (\tau)
      -
      \nabla\left(\div c_{i+1} (\tau)\right)
    }_{\L1 (\reali^n; \reali^n)}
    \d\tau
    \Bigr]
    \\
    & &
    +
    e^{\kappa^* t} \norma{u_o}_{\L\infty (\reali^n;\reali)}
    \\
    & &
    \qquad
    \times
    \int_{0}^t
    e^{(\alpha+K) \norma{w_o}_{\L\infty (\reali^n; \reali)}\tau e^{\gamma \tau}}
    \norma{
      b_{i+1} (\tau) - b_i (\tau)
      -
      \div\left(c_{i+1} (\tau) - c_i (\tau)\right)
    }_{\L{1} (\reali^n; \reali)}
    \d\tau .
  \end{eqnarray*}
  We proceed estimating all terms appearing in the inequality above.
  Begin by 5.~in Proposition~\ref{prop:hyper}, {\vv} and Claim~1:
  \begin{eqnarray}
    \kappa^*
    & = &
    \norma{b_{i+1}}_{\L\infty ([0,t] \times \reali^n;\reali)}
    + \norma{\div (c_{i+1} - c_i )}_{\L\infty ([0,t] \times \reali^n;\reali)}
    \nonumber
    \\
    & \leq &
    \alpha  \, \norma{w_o}_{\L\infty (\reali^n; \reali)} \, e^{\gamma t} + \beta
    +
    \norma{\div c_{i+1}}_{\L\infty ([0,t] \times \reali^n;\reali)}
    +
    \norma{\div c_i}_{\L\infty ([0,t] \times \reali^n;\reali)}
    \nonumber
    \\
    & \leq &
    \alpha  \, \norma{w_o}_{\L\infty (\reali^n; \reali)} \, e^{\gamma t}+ \beta
    +
    2 K \, \norma{w_o}_{\L\infty (\reali^n; \reali)}  \, e^{\gamma t}
    \nonumber
    \\
    & = &
    \lambda
    \qquad \mbox{ where } \qquad
    \lambda
    =
    (\alpha + 2 K)  \, \norma{w_o}_{\L\infty (\reali^n; \reali)} \, e^{\gamma t}
    + \beta
    \,;
    \label{eq:kappa}
    \\
    \nonumber
    \kappa_1^*
    & = &
    \norma{b_{i+1}}_{\L\infty ([0,t] \times \reali^n;\reali)}
    + \norma{\nabla \cdot (c_{i+1} - c_i )}_{\L\infty ([0,t] \times \reali^n;\reali)}
    + (2n+1) \norma{\nabla c_{i+1}}_{\L\infty ([0,t] \times \reali^n; \reali)}
    \\
    & \leq &
    (\alpha + 2 K) \, \norma{w_o}_{\L\infty (\reali^n; \reali)}  \, e^{\gamma t}
    + \beta
    +
    (2n+1) \, K \,  \norma{w_o}_{\L\infty (\reali^n; \reali)} \, e^{\gamma t}
    \nonumber
    \\
    & = &
    \lambda_1
    \qquad \mbox{ where } \qquad
    \lambda_1
    =
    \left( \alpha + (2n + 3) \, K \right)  \,
    \norma{w_o}_{\L\infty (\reali^n; \reali)} \, e^{\gamma t}
    +
    \beta
    \,.
    \label{eq:kappa1}
  \end{eqnarray}
  By~{\vv} and Claim~1
  \begin{eqnarray}
    \nonumber
    \norma{c_{i+1}-c_i}_{\L1 ([0,t]; \L\infty (\reali^n; \reali^n))}
    & = &
    \norma{v (w_i) - v (w_{i-1})}_{\L1 ([0,t]; \L\infty (\reali^n; \reali^n))}
    \\
    \label{eq:diffcn}
    & \leq &
    K \, \norma{w_i - w_{i-1}}_{\L1 ([0,t]; \L1 (\reali^n;\reali))} .
  \end{eqnarray}
  Consider now
  \begin{eqnarray}
    \nonumber
    \!\!\!
    & &
    \!\!\!
    \int_0^t
    e^{(\alpha + K) \, \norma{w_o}_{\L\infty (\reali^n; \reali)} \, \tau \, e^{\gamma \tau}}
    \, \norma{\nabla b_{i+1} (\tau)
      - \nabla \left(\div c_{i+1} (\tau)\right)}_{\L1 (\reali^n; \reali^n)} \,
    \d\tau
    \\
    \!\!\!
    & \leq  &
    \!\!\!
    e^{(\alpha + K) \, \norma{w_o}_{\L\infty (\reali^n; \reali)} \, t \, e^{\gamma t}}
    \, \norma{\nabla b_{i+1}
      - \nabla \left(\div c_{i+1}\right)}_{\L1 ([0,t] \times \reali^n;\reali^n)}
    \nonumber
    \\
    \label{eq:int}
    \!\!\!
    & \leq  &
    \!\!\!
    e^{(\alpha + K) \, \norma{w_o}_{\L\infty (\reali^n; \reali)} \, t \, e^{\gamma t}} \!
    \left[
      \norma{\nabla b_{i+1}}_{\L1 ([0,t] \times \reali^n;\reali^n)}
      + \norma{\nabla \left(\div c_{i+1}\right)}_{\L1 ([0,t] \times \reali^n;\reali^n)}
    \right]
  \end{eqnarray}
  and observe that $\norma{\nabla b_{i+1}}_{\L1 ([0,t] \times
    \reali^n; \reali^n)} = \alpha \, \norma{\nabla w_i}_{\L1 ([0,t]
    \times \reali^n; \reali^n)} $. Recall the proof of~9.~in
  Proposition~\ref{prop:para}, use Claim~1, \eqref{eq:an} and the
  expression~\eqref{eq:3P} of $J_n$:
  \begin{eqnarray*}
    & &\norma{\nabla w_i (\tau)}_{\L1 (\reali^n; \reali^n)}
    \\
    & \leq &
    \norma{ \nabla H_\mu (\tau)}_{\L1 (\reali^n; \reali)}
    \norma{w_o}_{\L1 (\reali^n; \reali)}
    \\
    & &
    + \int_0^{\tau} \norma{ \nabla H_\mu (\tau - s)}_{\L1
      (\reali^n;\reali)} \, \norma{a_i (s)}_{\L\infty (\reali^n;
      \reali)}\, \norma{w_i (s)}_{\L1 (\reali^n; \reali)} \d{s}
    \\
    & \leq &
    \norma{ \nabla H_\mu (\tau)}_{\L1 (\reali^n; \reali)}
    \norma{w_o}_{\L1 (\reali^n; \reali)}
    \\
    & &
    + \norma{a_i}_{\L\infty ([0,\tau] \times \reali^n; \reali)}\,
    \int_0^{\tau} \norma{ \nabla H_\mu (\tau - s)}_{\L1 (\reali^n; \reali)} \,
    e^{\gamma s} \, \norma{w_o}_{\L1 (\reali^n; \reali)}
    \,\d{s}
    \\
    & \leq &
    \norma{ \nabla H_\mu (\tau)}_{\L1 (\reali^n; \reali)}
    \norma{w_o}_{\L1 (\reali^n; \reali)} + \norma{w_o}_{\L1 (\reali^n; \reali)} \,
    e^{\gamma \tau} \,
    \int_0^{\tau} \norma{ \nabla H_\mu (\tau - s)}_{\L1 (\reali^n; \reali)} \, \d{s}
    \\
    & &
    \times
    \left( \gamma +
      \delta \, \norma{u_o}_{\L\infty(\reali^n;\reali)}
      e^{ (\alpha + K) \,
        \norma{w_o}_{\L\infty (\reali^n; \reali)} \tau \,  e^{\gamma \tau}
      }
    \right)
    \\
    & \leq  &
    \norma{w_o}_{\L1 (\reali^n; \reali)} \biggl[ \frac{J_n}{\sqrt{\mu \tau}}
    + e^{\gamma \tau}\!
    \left( \gamma +
      \delta \, \norma{u_o}_{\L\infty(\reali^n;\reali)}
      e^{(\alpha + K) \,
        \norma{w_o}_{\L\infty (\reali^n; \reali)} \tau \, e^{\gamma \tau}
      }
    \right) \!
    \int_0^\tau \!\!\frac{J_n}{\sqrt{\mu (\tau - s)}} \d{s} \biggr]
    \\
    & \leq &
    \norma{w_o}_{\L1 (\reali^n; \reali)} \frac{J_n}{\sqrt{\mu \tau}}
    \left[
      1
      +
      2 \tau e^{\gamma \tau}
      \left(
        \gamma
        +
        \delta \norma{u_o}_{\L\infty(\reali^n;\reali)}
        e^{
          (\alpha + K)
          \norma{w_o}_{\L\infty (\reali^n; \reali)} \tau e^{\gamma \tau}
        }
      \right)
    \right].
  \end{eqnarray*}
  Therefore
  \begin{equation}
    \label{eq:gradw}
    \begin{array}{@{}l@{\,\,}l}
      \norma{\nabla b_{i+1}}_{\L1 ([0,t] \times \reali^n; \reali^n)}
      \leq &
      \alpha \norma{w_o}_{\L1 (\reali^n; \reali)} J_n
      \dfrac{\sqrt{t}}{\sqrt{\mu}}
      \\
      &
      \times
      \left[ 1 + 2 t e^{\gamma t}
        \left( \gamma +
          \delta \norma{u_o}_{\L\infty(\reali^n; \reali)}
          e^{ (\alpha + K) \norma{w_o}_{\L\infty (\reali^n; \reali)}
            t e^{\gamma t}
          }
        \right)
      \right].
    \end{array}
  \end{equation}
  Consider $ \norma{ \nabla \left( \div c_{i+1} \right)}_{\L1 ([0,t]
    \times \reali^n; \reali^n)} = \norma{ \nabla \left( \div v (w_i)
    \right)}_{\L1 ([0,t] \times \reali^n; \reali^n)}$. By~{\vv} and
  Claim~1:
  \begin{eqnarray}
    \nonumber
    \norma{\nabla \left( \div c_{i+1} \right)}_{\L1 ([0,t] \times \reali^n;\reali^n)}
    & \leq &
    \norma{C \left(\norma{w_i}_{\L1 (\reali^n; \reali)}\right) \,
      \norma{w_i}_{\L1 (\reali^n; \reali)}}_{\L1 ([0,t]; \reali)}
    \\
    & \leq &
    \norma{C \left(\norma{w_i}_{\L1 (\reali^n; \reali)}\right) }_{\L\infty ([0,t]; \reali)} \,
    \norma{w_i}_{\L1 ([0,t] \times \reali^n; \reali)}
    \nonumber
    \\
    & \leq &
    C \left( \norma{w_i}_{\L\infty ([0,t];\L1 (\reali^n; \reali))}  \right) \,
    \norma{w_i}_{\L1 ([0,t] \times \reali^n; \reali)}
    \nonumber
    \\
    \label{eq:graddivc}
    & \leq &
    C \left( \norma{w_o}_{\L1 (\reali^n; \reali)} \, e^{\gamma t}  \right) \,
    \norma{w_o}_{\L1 (\reali^n; \reali)} \, t \, e^{\gamma t}.
  \end{eqnarray}
  Now use~\eqref{eq:gradw} and~\eqref{eq:graddivc} in~\eqref{eq:int}:
  \begin{eqnarray}
    \nonumber
    & &
    \int_0^t
    e^{(\alpha + K) \, \norma{w_o}_{\L\infty (\reali^n; \reali)} \, \tau \, e^{\gamma \tau}}
    \, \norma{\nabla b_{i+1} (\tau)
      - \nabla \left(\div c_{i+1} (\tau)\right)}_{\L1 (\reali^n; \reali^n)} \,
    \d\tau
    \\
    & \leq &
    e^{(\alpha + K) \, \norma{w_o}_{\L\infty (\reali^n; \reali)} \, t \, e^{\gamma t}} \,
    \biggl[ C \left( \norma{w_o}_{\L1 (\reali^n; \reali)} \, e^{\gamma t} \right) \,
    \norma{w_o}_{\L1 (\reali^n; \reali)} \, t \, e^{\gamma t}
    \nonumber
    \\
    & &
    +
    \alpha \, \norma{w_o}_{\L1 (\reali^n; \reali)} \, J_n \,
    \frac{\sqrt{t}}{\sqrt{\mu}}
    \left(
      1 + 2 \, t \, e^{\gamma t}
      \left(
        \gamma
        + \delta \, \norma{u_o}_{\L\infty(\reali^n; \reali)}
        e^{ (\alpha + K) \,
          \norma{w_o}_{\L\infty (\reali^n; \reali)} \,  t \, e^{\gamma t}
        }
      \right)
    \right)
    \biggr]
    \nonumber
    \\
    \label{eq:bleahOK}
    & = &
    \boldsymbol{M}_t (u_o, w_o),
  \end{eqnarray}
  where, for brevity, we set
  \begin{align}
    \nonumber & \boldsymbol{M}_t (u_o, w_o)
    \\
    \label{eq:M}
    \begin{split}
      = {}\,\, & \norma{w_o}_{\L1 (\reali^n; \reali)} \, e^{(\alpha +
        K) \, \norma{w_o}_{\L\infty (\reali^n; \reali)} t \, e^{\gamma
          t}} \sqrt{t} \biggl[ C \! \left( \norma{w_o}_{\L1 (\reali^n;
          \reali)} e^{\gamma t} \right) \sqrt{t} \, e^{\gamma t}
      \\
      & + J_n \, \dfrac{\alpha}{\sqrt{\mu}} \left( 1 + 2 t e^{\gamma
          t} \left( \gamma + \delta \norma{u_o}_{\L\infty(\reali^n;
            \reali)} e^{(\alpha + K) \, \norma{w_o}_{\L\infty
              (\reali^n; \reali)} t \, e^{\gamma t}} \right) \right)
      \biggr]
    \end{split}
    \\
    \label{eq:Mr}
    \leq \,{}\, & r\, e^{(\alpha + K) \, r \, t \, e^{\gamma t}}
    \sqrt{t} \, \biggl[ C \! \left( r \, e^{\gamma t} \right) \sqrt{t}
    \, e^{\gamma t} + J_n \, \dfrac{\alpha}{\sqrt{\mu}} \left( 1 + 2 t
      e^{\gamma t} \left( \gamma + \delta r \, e^{(\alpha + K) \, r \,
          t \, e^{\gamma t}} \right) \right) \biggr],
  \end{align}
  since $(u_o, w_o) \in \mathcal{X}_r$. Pass to
  \begin{eqnarray}
    \nonumber
    & &
    \int_{0}^t
    e^{(\alpha+K) \norma{w_o}_{\L\infty (\reali^n; \reali)}\tau e^{\gamma \tau}}
    \norma{
      b_{i+1} (\tau) - b_i (\tau)
      -
      \div\left(c_{i+1} (\tau) - c_i (\tau)\right)
    }_{\L{1} (\reali^n; \reali)}
    \d\tau
    \\
    \label{eq:uffa0}
    & \leq &
    e^{(\alpha+K) \norma{w_o}_{\L\infty (\reali^n; \reali)}t e^{\gamma t}}
    \!\!\!\int_0^t
    \norma{
      b_{i+1} (\tau) - b_i (\tau)
      -
      \div\left(c_{i+1} (\tau) - c_i (\tau)\right)
    }_{\L{1} (\reali^n; \reali)}
    \!\!\!\d\tau \,.
  \end{eqnarray}
  In particular, by the definition of $b_i$, {\vv} and Claim~1, we
  have
  \begin{eqnarray}
    \nonumber
    & &
    \int_0^t
    \norma{
      b_{i+1} (\tau) - b_i (\tau)
      -
      \div\left(c_{i+1} (\tau) - c_i (\tau)\right)
    }_{\L{1} (\reali^n; \reali)}
    \d\tau
    \\
    \nonumber
    & \leq &
    \int_0^t
    \norma{ b_{i+1} (\tau) - b_i (\tau)}_{\L1 (\reali^n; \reali)} \d\tau
    +
    \int_0^t
    \norma{ \div\left(c_{i+1} (\tau) - c_i (\tau)\right)}_{\L1 (\reali^n; \reali)}
    \d\tau
    \\
    \nonumber
    & \leq &
    \alpha
    \int_0^t \norma{ w_i (\tau) - w_{i-1} (\tau)}_{\L1 (\reali^n; \reali)} \d\tau
    \\
    \nonumber
    & &
    +
    \int_0^t
    C\left( \norma{w_{i-1} (\tau)}_{\L\infty (\reali^n; \reali)} \right)
    \norma{w_i (\tau) - w_{i-1} (\tau)}_{\L1 (\reali^n; \reali)} \, \d\tau
    \notag
    \\
    \nonumber
    & \leq &
    \alpha \, \norma{w_i - w_{i-1}}_{\L1 ([0,t];\L1 (\reali^n; \reali))}
    \\
    \nonumber
    & &
    +
    \norma{C\left( \norma{w_{i-1} (\tau)}_{\L\infty (\reali^n; \reali)} \right)}_{\L\infty ([0,t]; \reali)} \,
    \norma{w_i - w_{i-1}}_{\L1 ([0,t];\L1 (\reali^n; \reali))}
    \\
    \nonumber
    & \leq &
    \left(
      \alpha + C\left(
        \norma{w_{i-1}}_{\L\infty ([0,t] \times \reali^n; \reali)} \right)
    \right) \,
    \norma{w_i - w_{i-1}}_{\L1 ([0,t];\L1 (\reali^n; \reali))}
    \\
    \label{eq:uffa1}
    & \leq &
    \left(
      \alpha + C\left(
        \norma{w_o}_{\L\infty (\reali^n; \reali)} \, e^{\gamma t} \right)
    \right)
    \norma{w_i - w_{i-1}}_{\L1 ([0,t];\L1 (\reali^n; \reali))} \, .
  \end{eqnarray}
  Use~\eqref{eq:kappa}, \eqref{eq:kappa1}, \eqref{eq:diffcn},
  \eqref{eq:bleahOK}, \eqref{eq:uffa0} and~\eqref{eq:uffa1} to bound
  $\norma{u_{i+1} (t) - u_i (t)}_{\L1 (\reali^n; \reali)}$:
  \begin{eqnarray*}
    & &
    \norma{u_{i+1} (t) - u_i (t)}_{\L1 (\reali^n; \reali)}
    \\
    & \leq &
    e^{\lambda \, t} \,
    \norma{w_i - w_{i-1}}_{\L1 ([0,t]; \L1 (\reali^n; \reali))}
    \biggl[
    K \, e^{\lambda_1 \, t} \,
    \left[
      \tv (u_o)
      +
      I_n \,  \norma{u_o}_{\L\infty (\reali^n; \reali)} \,
      \boldsymbol{M}_t (u_o, w_o)
    \right]
    \\
    & &
    +
    \norma{u_o}_{\L\infty (\reali^n; \reali)} \,
    \left(
      \alpha + C\left(
        \norma{w_o}_{\L\infty (\reali^n; \reali)} \, e^{\gamma t} \right)
    \right)
    e^{
      (\alpha + K) \,
      \norma{w_o}_{\L\infty (\reali^n; \reali)} \, t \, e^{\gamma t}
    }
    \biggr]
  \end{eqnarray*}
  and then, recalling~\eqref{eq:kappa},~\eqref{eq:kappa1}
  and~\eqref{eq:M},
  \begin{eqnarray}
    \nonumber
    & &
    \norma{u_{i+1} - u_i }_{\L1 ([0,T]; \L1(\reali^,; \reali))}
    \\
    \nonumber
    & \leq &
    T \,
    \norma{w_i - w_{i-1}}_{\L1 ([0,T]; \L1 (\reali^n;\reali))}
    \exp
    \left( (\alpha + 2 K) \,
      \norma{w_o}_{\L\infty (\reali^n; \reali)} \, T \, e^{\gamma T}
      + \beta \, T
    \right)
    \\
    \nonumber
    & &
    \times
    \Biggl\{
    K \, \exp
    \left( \left( \alpha + (2n + 3) \, K  \right)
      \norma{w_o}_{\L\infty (\reali^n; \reali)} \, T \,
      e^{\gamma T} + \beta \, T
    \right)
    \biggl\{ \tv (u_o)
    \\
    \nonumber
    & &
    \quad
    +
    I_n \, \norma{u_o}_{\L\infty (\reali^n; \reali)} \,
    \norma{w_o}_{\L1 (\reali^n; \reali)} \, \sqrt{T}
    e^{
      (\alpha + K) \, \norma{w_o}_{\L\infty (\reali^n; \reali)} \,
      T \, e^{\gamma T}
    }
    \\
    \nonumber
    & &
    \quad
    \times
    \biggl[ J_n \frac{\alpha}{\sqrt{\mu}}
    \left(
      1 + 2 \, T \, e^{\gamma T}\!
      \left(
        \gamma + \delta \, \norma{u_o}_{\L\infty(\reali^n; \reali)}
        e^{
          (\alpha + K) \, \norma{w_o}_{\L\infty (\reali^n; \reali)} \, T \,
          e^{\gamma T}
        }
      \right)
    \right)
    \\
    \nonumber
    & &
    \qquad
    + C \left( \norma{w_o}_{\L1 (\reali^n; \reali)} \,
      e^{\gamma T} \right) \, \sqrt{T} \, e^{\gamma T}
    \biggr]
    \biggr\}
    \\
    \nonumber
    & &
    \quad
    +
    \norma{u_o}_{\L\infty (\reali^n; \reali)}
    \left(
      \alpha
      +
      C\left( \norma{w_o}_{\L\infty (\reali^n; \reali)} \, e^{\gamma T}
      \right)\!
    \right)
    e^{
      (\alpha + K) \, \norma{w_o}_{\L\infty (\reali^n; \reali)} T \,
      e^{\gamma T} }
    \Biggr\}
    \\
    \label{eq:diffuOK}
    & \leq &
    T \, \mathcal{K}_u (T,r) \,
    \norma{w_i - w_{i-1}}_{\L1 ([0,T]; \L1 (\reali^n;\reali))},
  \end{eqnarray}
  where, by~\eqref{eq:Mr},
  \begin{align}
    \label{eq:Kut}
    \begin{split}
      \mathcal{K}_u (T,r) \, = {}\, & r \, e^{\left((2 \alpha + 3 K)
          \, r \, e^{\gamma T} + \beta\right) T} \Biggl\{ K \,
      e^{\left(2(n + 1) \, K \, r \, e^{\gamma T} + \beta\right) T}
      \\
      & \times \biggl\{ 1 + I_n \, r \sqrt{T} e^{(\alpha + K) r T \,
        e^{\gamma T}} \biggl[ J_n \frac{\alpha}{\sqrt{\mu}} \left( 1 +
        2 T e^{\gamma T}\!  \left( \gamma + \delta \, r \, e^{ (\alpha
            + K) \, r \, T \, e^{\gamma T} } \right) \right)
      \\
      & \qquad + C \left( r \, e^{\gamma T} \right) \, \sqrt{T} \,
      e^{\gamma T} \biggr] \biggr\} + \alpha + C\left( r \, e^{\gamma
          T} \right)\!  \Biggr\}.
    \end{split}
  \end{align}
  Therefore, the above definition of $\mathcal{K}_u (T,r)$, together
  with~\eqref{eq:diffw}, \eqref{eq:KTW} and~\eqref{eq:diffuOK}, yields
  \begin{eqnarray*}
    & &
    d\left( \left(u_{i+1}, w_{i+1}\right), \left(u_i, w_i\right)
    \right)
    \\
    & =\! &
    \norma{u_{i+1} - u_i}_{\L1 ([0,T]; \L1 (\reali^n; \reali))} +
    \norma{w_{i+1} - w_i}_{\L1 ([0,T]; \L1 (\reali^n; \reali))}
    \\
    & \leq\! &
    T
    \max
    \left\{ \mathcal{K}_u (T,r) , \, \mathcal{K}_w (T,r)\right\}\!
    \left[
      \norma{u_i - u_{i-1}}_{\L1 ([0,T]; \L1 (\reali^n;\reali))}
      +
      \norma{w_i - w_{i-1}}_{\L1 ([0,T]; \L1(\reali^n;\reali))}
    \right]
    \\
    & =\! &
    T \,
    \mathcal{K} (T,r) \,
    d\left( \left(u_i, w_i\right) , \left(u_{i-1}, w_{i-1} \right)\right)
  \end{eqnarray*}
  where
  \begin{equation}
    \label{eq:K}
    \mathcal{K} (T,r)
    =
    \max \left\{
      \mathcal{K}_u (T,r), \,
      \mathcal{K}_w (T,r)
    \right\} ,
  \end{equation}
  proving~\eqref{eq:distanza}.

  For any positive $r$, we can now choose $T_r$ so that $T_r \,
  \mathcal{K} (T_r, r) < 1$. The sequence $(u_i, w_i)$ converges in $X
  = \L1 ([0, T_r]; \mathcal{X}^+)$ to a limit, say, $(u_*, w_*)$. By
  construction, see Claim~0, both $u_*$ and $w_*$ attain non negative
  values. We now check that $(u_*, w_*)$ solves~\eqref{eq:Model} in
  the sense of Definition~\ref{def:Main}.

  Clearly, $(u_*, w_*) (0) = (u_o, w_o)$. Moreover, by the above
  construction, we have that for any $\phi \in \Cc\infty (\left]0, T_r
  \right[ \times \reali^n; \reali)$
  \begin{eqnarray*}
    \int_0^{T_r} \int_{\reali^n}
    \left(
      w_i \, \partial_t \phi
      +
      \mu \, w_i \, \Delta\phi
      +
      (\gamma - \delta \, u_{i-1}) \, w_i \, \phi
    \right)
    \d{x} \d{t}
    & \!\!=\!\! &
    0
    \\
    \int_0^{T_r} \int_{\reali^n}
    \left(
      u_i \, \partial_t \phi
      +
      u_i \, v (w_{i-1}) \cdot \nabla\phi
      +
      (\alpha\, w_{i-1} - \beta) \, u_i \, \phi
    \right)
    \d{x} \d{t}
    & \!\!=\!\! &
    0 \,.
  \end{eqnarray*}
  Thanks to the $\L\infty$ bounds proved in Claims~1 and~2, we can
  apply the Dominated Convergence Theorem, ensuring that $(u_*, w_*)$
  is a weak solution to~\eqref{eq:Model} for $t \in [0, T_r]$.

  For all $t \in [0,T_r]$, we define $\mathcal{R}_{0,t} (u_o,w_o) =
  (u_*, w_*) (t)$.

  \bigskip

  Consider now a couple of initial data $(u_{1,o},w_{1,o}), \,
  (u_{2,o},w_{2,o})$ satisfying~\eqref{eq:woRegular}. For $t \in
  [0,T_r]$, we know that $\left(u_i, w_i \right) (t) =
  \mathcal{R}_{0,t} (u_{i,o}, w_{i,o})$ solves~\eqref{eq:Model} with
  initial datum $(u_{i,o}, w_{i,o})$ in distributional sense. For
  $(t,x) \in [0,T_r] \times \reali^n$ we define for $i=1,2$
  \begin{equation*}
    a_i (t,x) = \gamma - \delta \, u_i (t,x)
    \,,\quad
    b_i (t,x) = \alpha\, w_i (t,x) - \beta
    \quad \mbox{ and } \quad
    c_i (t,x) = v\left(w_i (t)\right) (x).
  \end{equation*}
  Using the operators $\mathcal{P}$ of Proposition~\ref{prop:para} and
  $\mathcal{H}$ of Proposition~\ref{prop:hyper}, observe that $w_i (t)
  = \mathcal{P}_{0,t}^{\,i}w_{i,o}$ and $u_i (t) =
  \mathcal{H}_{0,t}^{\,i}u_{i,o}$, for $i=1,2$. Moreover, note that
  $\mathcal{P}_{0,t}^1w_{2,o}$ is the solution to~\eqref{eq:CauchyP}
  with $a_1$ in the source term and initial datum $w_{2,o}$, while
  $\mathcal{H}_{0,t}^1u_{2,o}$ is the solution to~\eqref{eq:CauchyH}
  with coefficients $b_1, \, c_1$ and initial datum $u_{2,o}$. We
  compute $ \norma{ \mathcal{R}_{0,t} (u_{1,o},w_{1,o}) -
    \mathcal{R}_{0,t} (u_{2,o},w_{2,o})}_{\mathcal{X}} $ as defined
  in~\eqref{eq:normaX}:
  \begin{eqnarray}
    \nonumber
    & &
    \norma{
      \mathcal{R}_{0,t} (u_{1,o},w_{1,o}) -
      \mathcal{R}_{0,t} (u_{2,o},w_{2,o})}_{\mathcal{X}}
    \\
    \nonumber
    & = &
    \norma{
      \mathcal{H}_{0,t}^1 u_{1,o} - \mathcal{H}_{0,t}^2 u_{2,o}
    }_{\L1 (\reali^n; \reali)}
    +
    \norma{
      \mathcal{P}_{0,t}^1 w_{1,o} - \mathcal{P}_{0,t}^2 w_{2,o}
    }_{\L1 (\reali^n; \reali)}
    \\
    & \leq &
    \label{eq:stimaRH}
    \norma{
      \mathcal{H}_{0,t}^1 u_{1,o} - \mathcal{H}_{0,t}^1 u_{2,o}
    }_{\L1 (\reali^n; \reali)}
    +
    \norma{
      \mathcal{H}_{0,t}^1 u_{2,o} - \mathcal{H}_{0,t}^2 u_{2,o}
    }_{\L1 (\reali^n; \reali)}
    \\
    \label{eq:stimaRP}
    & &
    +
    \norma{
      \mathcal{P}_{0,t}^1 w_{1,o} - \mathcal{P}_{0,t}^1 w_{2,o}
    }_{\L1 (\reali^n; \reali)}
    +
    \norma{
      \mathcal{P}_{0,t}^1 w_{2,o} - \mathcal{P}_{0,t}^2 w_{2,o}
    }_{\L1 (\reali^n; \reali)} .
  \end{eqnarray}
  Compute each term of~\eqref{eq:stimaRH} separately.  Since the map
  $\mathcal{H}_{0,t}^1$ is linear, by~3.~in
  Proposition~\ref{prop:hyper} and its particularization in Claim~2,
  the first term in~\eqref{eq:stimaRH} can be estimate as follows:
  \begin{eqnarray}
    \nonumber
    \norma{
      \mathcal{H}_{0,t}^1 u_{1,o}
      -
      \mathcal{H}_{0,t}^1 u_{2,o} }_{\L1 (\reali^n; \reali)}
    & \leq &
    \norma{u_{1,o} - u_{2,o}}_{\L1 (\reali^n; \reali)}
    \exp
    \left(
      \alpha \,
      \norma{w_{1,o}}_{\L\infty (\reali^n;\reali)} \,
      t \,
      e^{\gamma \, t}
    \right)
    \\
    \label{eq:h1}
    & \leq &
    \norma{u_{1,o} - u_{2,o}}_{\L1 (\reali^n; \reali)}
    \exp \left( \alpha \, r \, t \, e^{\gamma \, t}  \right).
  \end{eqnarray}
  Concerning the second term in~\eqref{eq:stimaRH}, recall~5.~in
  Proposition~\ref{prop:hyper} and adapt the estimates above for
  $\norma{u_{i+1} (t) - u_i (t)}_{\L1 (\reali^n; \reali)}$, using
  $\boldsymbol{M}_t$ as defined in~\eqref{eq:M}--\eqref{eq:Mr} and
  \begin{equation*}
    \begin{array}{l@{}l}
      \Theta
      & =
      \left(
        (\alpha + K) \, \norma{w_{1,o}}_{\L\infty (\reali^n; \reali)}
        + K \norma{w_{2,o}}_{\L\infty (\reali^n; \reali)}
      \right)
      e^{\gamma \, t}
      + \beta ,
      \\
      \Theta_1
      & =
      \left(
        (\alpha + 2 \, (n+1) \, K) \, \norma{w_{1,o}}_{\L\infty (\reali^n; \reali)}
        + K \norma{w_{2,o}}_{\L\infty (\reali^n; \reali)}
      \right)
      e^{\gamma \, t}
      + \beta.
    \end{array}
  \end{equation*}
  obtaining
  \begin{eqnarray}
    \nonumber
    & &
    \norma{
      \mathcal{H}_{0,t}^1 u_{2,o} - \mathcal{H}_{0,t}^2 u_{2,o}
    }_{\L1 (\reali^n; \reali)}
    \\
    & \leq &
    e^{\Theta \, t} \,
    \int_0^t
    \norma{\mathcal{P}_{0,\tau}^1w_{1,o} -
      \mathcal{P}_{0,\tau}^2w_{2,o}}_{\L1 (\reali^n; \reali)} \, \d\tau
    \nonumber
    \\
    \nonumber
    & &
    \times
    \biggl[
    K \, e^{\Theta_1 \, t} \,
    \left[
      \tv (u_{2,o})
      +
      I_n \, \norma{u_{2,o}}_{\L\infty (\reali^n; \reali)} \,
      \boldsymbol{M}_t (u_{1,o}, w_{1,o})
    \right]
    \\
    & &
    + \norma{u_{2,o}}_{\L\infty (\reali^n; \reali)} \,
    \left(
      \alpha + C \left(
        \norma{w_{2,o}}_{\L\infty (\reali^n; \reali)} \, e^{\gamma \, t}
      \right)
    \right)
    e^{
      (\alpha + K) \, t \, e^{\gamma \, t} \,
      \underset{i=1,2}{\max} \norma{w_{i,o}}_{\L\infty (\reali^n; \reali)}
    }
    \biggr]
    \nonumber
    \\
    \label{eq:h2}
    & \leq &
    \mathcal{K}_u (t,r)
    \int_0^t
    \norma{\mathcal{P}_{0,\tau}^1w_{1,o} -
      \mathcal{P}_{0,\tau}^2w_{2,o}}_{\L1 (\reali^n; \reali)} \, \d\tau \,,
  \end{eqnarray}
  where $\mathcal{K}_u (t,r)$ is defined in~\eqref{eq:Kut}.

  Pass to~\eqref{eq:stimaRP}.  Since the map $\mathcal{P}_{0,t}^1$ is
  linear, by~3.~in Proposition~\ref{prop:para} and its
  particularization in Claim~1, we have the following estimate for the
  first term in~\eqref{eq:stimaRP}:
  \begin{equation}
    \label{eq:p1}
    \norma{
      \mathcal{P}_{0,t}^1 w_{,o} - \mathcal{P}_{0,t}^1 w_{2,o}
    }_{\L1 (\reali^n; \reali)}
    \leq
    \norma{w_{1,o} - w_{2,o}}_{\L1 (\reali^n; \reali)} \,
    e^{\gamma \, t}.
  \end{equation}
  Concerning the second term in~\eqref{eq:stimaRP}, recall~5.~in
  Proposition~\ref{prop:para} and adapt the estimates above for
  $\norma{w_{i+1} (t) - w_i (t)}_{\L1 (\reali^n; \reali)}$ to obtain
  \begin{eqnarray}
    \nonumber
    & &
    \norma{ \mathcal{P}_{0,t}^1 w_{2,o} - \mathcal{P}_{0,t}^2
      w_{2,o} }_{\L1 (\reali^n; \reali)}
    \\
    \nonumber
    & = &
    \delta \, \norma{w_{2,o}}_{\L\infty (\reali^n; \reali)} \,
    \exp \left( 2 \,\gamma \, t + \delta \, t \,
      \norma{u_{2,o}}_{\L\infty (\reali^n; \reali)} \,
      e^{
        (\alpha + K) \, \norma{w_{2.o}}_{\L\infty (\reali^n; \reali)} \,
        \frac{e^{\gamma \,t}-1}{\gamma}
      }
    \right)
    \\
    \nonumber
    & &
    \times \int_0^t \norma{\mathcal{H}_{0,\tau}^1 u_{1,o} -
      \mathcal{H}_{0,\tau}^2 u_{2,o}}_{\L1 (\reali^n; \reali)} \,
    \d\tau
    \\
    \label{eq:p2}
    & \leq &
    \mathcal{K}_w (t,r)
    \int_0^t \norma{\mathcal{H}_{0,\tau}^1 u_{1,o} -
      \mathcal{H}_{0,\tau}^2 u_{2,o}}_{\L1 (\reali^n; \reali)}
    \d\tau
  \end{eqnarray}
  with $\mathcal{K}_w (t,r)$ as in~\eqref{eq:KTW}.  We can now
  rewrite~\eqref{eq:stimaRH}--\eqref{eq:stimaRP}
  using~\eqref{eq:h1}--\eqref{eq:h2}--\eqref{eq:p1}--\eqref{eq:p2}:
  \begin{eqnarray*}
    & &
    \norma{
      \mathcal{R}_{0,t} (u_{1,o},w_{1,o}) -
      \mathcal{R}_{0,t} (u_{2,o},w_{2,o})}_{\mathcal{X}}
    \\
    & = &
    \exp\left(\alpha \, r \, t \, e^{\gamma \, t} \right) \,
    \norma{u_{1,o} - u_{2,o}}_{\L1 (\reali^n; \reali)}
    +
    \mathcal{K}_u (t,r) \int_0^t
    \norma{
      \mathcal{P}_{0,\tau}^1 w_{1,o}
      -
      \mathcal{P}_{0,t}^2 w_{2,o}}_{\L1 (\reali^n; \reali)} \,
    \d\tau
    \\
    & &
    +
    e^{\gamma \, t} \,
    \norma{w_{1,o} - w_{2,o}}_{\L1 (\reali^n; \reali)}
    +
    \mathcal{K}_w (t,r)
    \int_0^t
    \norma{\mathcal{H}_{0,\tau}^1 u_{1,o} - \mathcal{H}_{0,t}^2 u_{2,o}
    }_{\L1 (\reali^n; \reali)} \, \d\tau
    \\
    & \leq &
    e^{(\gamma + \alpha \, r \, e^{\gamma \, t})\,t}
    \left(
      \norma{u_{1,o} - u_{2,o}}_{\L1 (\reali^n; \reali)}
      +
      \norma{w_{1,o} - w_{2,o}}_{\L1 (\reali^n; \reali)}
    \right)
    \\
    & &
    +
    \mathcal{K} (t,r)
    \int_0^t
    \norma{
      \mathcal{R}_{0,\tau} (u_{1,o},w_{1,o}) -
      \mathcal{R}_{0,\tau} (u_{2,o},w_{2,o})}_{\mathcal{X}} \, \d\tau
  \end{eqnarray*}
  and $\mathcal{K} (t,r)$ is as in~\eqref{eq:K}.  An application of
  Gronwall Lemma yields: for all $t\in [0, T_r]$
  \begin{equation}
    \label{eq:Lip}
    \begin{array}{c}
      \displaystyle
      \norma{
        \mathcal{R}_{0,t} (u_{1,o},w_{1,o}) -
        \mathcal{R}_{0,t} (u_{2,o},w_{2,o})}_{\mathcal{X}}
      \leq
      \mathcal{L} (t,r) \,
      \norma{(u_{1,o}, w_{1,o}) - (u_{2,o}, w_{2,o})}_{\mathcal{X}} \,,
      \\[3pt]
      \displaystyle
      \mbox{where} \qquad
      \mathcal{L} (t,r)
      =
      \exp
      \left[
        \left(
          \gamma
          +
          \alpha \, r \, e^{\gamma \, t}
          +
          \mathcal{K} (t,r)
        \right)t
      \right]
      \quad \mbox{ with } \quad
      \mathcal{K} (t,r) \mbox{ as in~\eqref{eq:K}} \,,
    \end{array}
  \end{equation}
  proving point~\ref{it:thm:Lipschitz}.~under
  condition~\eqref{eq:woRegular}.

  Denote $\left(u (t), w (t)\right) = \mathcal{R}_{0,t} (u_o,w_o)$ and
  define
  \begin{equation}
    \label{eq:LaSolita}
    a (t,x) = \gamma - \delta u (t,x)
    \,,\quad
    b (t,x) = \alpha\, w (t,x) - \beta
    \quad \mbox{ and } \quad
    c (t,x) = v\left(w (t)\right) (x).
  \end{equation}
  Then, with the notation in Proposition~\ref{prop:para} and
  Proposition~\ref{prop:hyper}, by construction
  \begin{equation}
    \label{eq:SempreLei}
    u (t) = \mathcal{H}_{0,t} u_o
    \quad \mbox{ and } \quad
    w (t) = \mathcal{P}_{0,t} w_o\,,
    \quad \mbox{ so that } \quad
    \mathcal{R}_{0,t} (u_o, w_o)
    =
    \left(
      \mathcal{P}_{0,t}w_o, \,
      \mathcal{H}_{0,t}u_o
    \right) \,.
  \end{equation}

  We now extend the map $t \to \mathcal{R}_{0,t} (u_o, w_o)$ to the
  whole time axis. Indeed, define the time
  \begin{displaymath}
    T_*
    =
    \sup
    \left\{
      T \in \reali^+ \colon
      \mathcal{R}_{0,t} (u_o, w_o) \mbox{ is defined for all } t \in [0,T]
    \right\} \,.
  \end{displaymath}
  The previous construction ensures that $\mathcal{R}_{0,t} (u_o,w_o)$
  is defined at least for all $t \in [0, T_r]$, hence the supremum
  above is well defined. Assume now that $T_* < +\infty$.  By~8.~in
  Proposition~\ref{prop:hyper}, the first component of the map $t \to
  \mathcal{R}_{0,t} (u_o, w_o)$ is Lipschitz continuous with a
  constant uniformly bounded by a function of $r$ on any bounded time
  interval.  By~8.~in Proposition~\ref{prop:para} the second component
  of the map $t \to \mathcal{R}_{0,t} (u_o, w_o)$ is uniformly
  continuous on, say $[0, T_r/2]$ and H\"older continuous with
  exponent $\theta=1/2$ for $t \in \left[T_r/2, T_*\right[$, the
  H\"older constant being bounded by a function of $r$.  Hence, the
  map $t \to \mathcal{R}_{0,t} (u_o, w_o)$ is uniformly continuous for
  $t \in \left[0, T_* \right[$ and the pair $(U,W) = \lim_{t\to T_*^-}
  \left(u (t), w (t)\right)$ is well defined and in
  $\mathcal{X}_{r_*}$ for a suitable $r_*$. Moreover, since $W = w
  (T_*) = \mathcal{P}_{0,T_*} w_o$, we have that $(U,W)$ also
  satisfies~\eqref{eq:woRegular} by 7.~and 9.~in
  Proposition~\ref{prop:para}. Repeating the construction above, we
  show that the Cauchy problem consisting of~\eqref{eq:Model} with
  initial datum $(u,w) (T_*) = (U,W)$ admits a solution defined on
  $[T_*, T_* + T_{r_*}]$, which contradicts the choice of $T_*$,
  unless $T_* = +\infty$.  Define $a,b$ and $c$ as
  in~\eqref{eq:LaSolita}. Then, by 8.~in Proposition~\ref{prop:para}
  and 8.~in Proposition~\ref{prop:hyper}, \eqref{eq:SempreLei}
  directly ensures the continuity in time of $\mathcal{R}$, for all $t
  \in \reali^+$. This completes the proof of
  point~\ref{it:thm:solution}.~under condition~\eqref{eq:woRegular}.

  Choose now a general initial datum $(u_o, w_o) \in \mathcal{X}^+$,
  so that $(u_o, w_o) \in \mathcal{X}_r$ for a suitable $r>0$. Let
  $\rho_n$ be a sequence of mollifiers with $\rho_n \in \Cc\infty
  (\reali^n; \reali^+)$ and $\int_{\reali^n} \rho_n (x) \d{x} =
  1$. Then, the sequence of initial data $(u_o, w_o * \rho_n)$ is in
  $\mathcal{X}_r$, satisfies~\eqref{eq:woRegular} and converges to
  $(u_o, w_o)$ in $\mathcal{X}$. By~\eqref{eq:Lip}, we can uniquely
  extend $\mathcal{R}$ through the limit $\mathcal{R}_{0,t} (u_o, w_o)
  = \lim_{n\to+\infty} \mathcal{R}_{0,t} (u_o, w_o * \rho_n)$ to all
  $\mathcal{X}_r$ and for all $t \in \reali^+$, completing the proof
  of point~\ref{it:thm:solution}.~and of
  point~\ref{it:thm:Lipschitz}.~for all $t \in \reali^+$. Note that
  the positivity of the solution follows from Claim~0. The $\L1$ and
  $\L\infty$ estimates at point~\ref{it:thm:growth}.~now follow from
  Claims~1 and~2. Again, the continuity in time of the map
  $\mathcal{R}$ so extended directly follows from 8.~in
  Proposition~\ref{prop:para} and 8.~in Proposition~\ref{prop:hyper}.

  We now prove that $\mathcal{R}$ is a process. For any $(u_o, w_o)
  \in \mathcal{X}^+$, use the notation~\eqref{eq:LaSolita} and observe
  that 1.~in Proposition~\ref{prop:para}, 1.~in
  Proposition~\ref{prop:hyper} and~\eqref{eq:SempreLei} ensure that
  the map $\mathcal{R}$ is a process. Note however
  that~\eqref{eq:Model} is autonomous, hence
  definitions~\ref{def:para}, \ref{def:Hyp} and~\ref{def:Main} ensure
  that $\mathcal{R}$ is a semigroup, proving point~\ref{it:thm:semi}.

  Point~\ref{thm:support}.~follows from~9.~in
  Proposition~\ref{prop:hyper}, {\vv} and the $\L1$ estimate in
  point~\ref{it:thm:growth}.~above.
\end{proofof}


\begin{lemma}
  \label{lem:v}
  Let $\eta$ be such that $\nabla\eta \in (\W{2}{1} \cap
  \W{1}{\infty}) (\reali^n; \reali^n)$. Then, the map $v$ defined
  in~\eqref{eq:w} satisfies~{\vv} with
  \begin{eqnarray}
    \label{eq:4}
    & &
    K
    =
    \max
    \left\{
      2 \kappa \, \norma{\nabla\eta}_{\W21},\,
      2 \kappa \, \norma{\nabla\eta}_{\W1\infty},\,
      3 \norma{\nabla\eta}_{\W1\infty},\,
      \frac{48}{25\sqrt{5}} \norma{\nabla\eta}_{\W21}
    \right\} ,
    \\
    \nonumber
    & &
    C\left( \xi \right)
    =
    K \,
    \left( 1 + K \, \xi \right)  \quad \mbox{for all } \xi \in \reali^+.
  \end{eqnarray}
\end{lemma}

\noindent In the proof below we use the Euclidean norm
$\norma{v}_{\reali^n} = \sqrt{\sum_{i=1}^n (v_i)^2}$ on vectors in
$\reali^n$ and the operator norm $\norma{A}_{\reali^{n\times n}} =
\sup_{v\colon \norma{v}_{\reali^n}=1} \norma{A\, v}_{\reali^n}$ on
$n\times n$ matrices.

\begin{proofof}{Lemma~\ref{lem:v}}
  The $\C2$ regularity is immediate. The bound on $\norma{v
    (w)}_{\L\infty (\reali^n; \reali)}$ is ensured by
  \begin{equation}
    \label{eq:v}
    \norma{v (w)}_{\L\infty (\reali^n; \reali^n)}
    \leq
    \kappa \, \norma{\nabla \eta}_{\L\infty (\reali^n; \reali^n)}
    \norma{w}_{\L1 (\reali^n; \reali)} ,
  \end{equation}
  see also~\cite[Lemma~3.1]{ColomboGaravelloLecureux}.  To estimate
  the gradient of $v (w)$, use the standard identity $\nabla (f\, v) =
  f \, \nabla v + v \otimes \nabla f$:
  \begin{eqnarray*}
    \nabla v (w)
    & = &
    \kappa \,
    \frac{1}{\left( 1 + \norma{w * \nabla \eta }^2 \right)^{1/2}}
    \nabla \left(w * \nabla \eta \right)
    +
    \kappa \, \left(w * \nabla \eta \right)
    \otimes
    \nabla
    \frac{1}{\left( 1 + \norma{w * \nabla \eta }^2 \right)^{1/2}}
    \\
    & = &
    \kappa \, \frac{ w * \nabla^2 \eta}
    {\left( 1 + \norma{w * \nabla \eta }^2 \right)^{1/2}}
    -
    \kappa \, \left(w * \nabla \eta \right)
    \otimes
    \frac{\left( w * \nabla^2 \eta \right) \left( w * \nabla \eta \right)}
    {\left( 1 + \norma{w * \nabla \eta }^2 \right)^{3/2}} \,.
  \end{eqnarray*}
  Recall that $\norma{v_1 \otimes v_2}_{\reali^{n\times n}} \leq
  \norma{v_1}_{\reali^n} \, \norma{v_2}_{\reali^n}$ and $\norma{A\,
    v}_{\reali^n} \leq \norma{A}_{\reali^{n\times n}} \,
  \norma{v}_{\reali^n}$, hence
  \begin{eqnarray}
    \nonumber
    & &
    \norma{\nabla v (w)}_{\L\infty (\reali^n; \reali^{n\times n})}
    \\
    \nonumber
    & \leq &
    \kappa \,
    \frac{\norma{ w * \nabla^2 \eta}_{\L\infty (\reali^n; \reali^{n\times n})}}
    {\sqrt{1 + \norma{w * \nabla \eta }^2}}
    \\
    \nonumber
    & &
    +
    \kappa \,
    \norma{
      \frac{w * \nabla \eta}
      {\sqrt{1 + \norma{w * \nabla \eta }^2}}
      \otimes
      \frac{w * \nabla^2 \eta }
      {\sqrt{1 + \norma{w * \nabla \eta }^2}}
      \frac{ w * \nabla \eta}
      {\sqrt{1 + \norma{w * \nabla \eta }^2}}
    }_{\L\infty (\reali^n; \reali^{n\times n})}
    \\
    \nonumber
    & \leq &
    \kappa \,
    \norma{ w * \nabla^2 \eta}_{\L\infty (\reali^n; \reali^{n\times n})}
    \\
    \nonumber
    & &
    +
    \kappa \,
    \norma{
      \frac{w * \nabla \eta}
      {\sqrt{1 + \norma{w * \nabla \eta }^2}}
    }_{\L\infty (\reali^n; \reali^n)}
    \!\!\! \norma{w * \nabla^2 \eta}_{\L\infty (\reali^n; \reali^{n\times n})}
    \norma{
      \frac{ w * \nabla \eta}
      {\sqrt{1 + \norma{w * \nabla \eta }^2}}
    }_{\L\infty (\reali^n; \reali^n)}
    \\
    \nonumber
    & \leq &
    2 \, \kappa \,
    \norma{w * \nabla^2 \eta}_{\L\infty (\reali^n; \reali^{n\times n})}
    \\
    \label{eq:gradv}
    & \leq &
    2 \, \kappa \,
    \norma{\nabla^2\eta}_{\L1 (\reali^n; \reali^{n\times n})}
    \norma{w}_{\L\infty (\reali^n;\reali)} \,,
  \end{eqnarray}
  proving also that $v (w) \in \W{1}{\infty} (\reali^n; \reali^n)$.
  Pass now to
  \begin{eqnarray*}
    & &
    v (w_1) - v (w_2)
    \\
    & = &
    \kappa \, \frac{(w_1-w_2)*\nabla\eta}{\sqrt{1+\norma{w_1*\nabla\eta}^2}}
    +
    \kappa \, (w_2*\nabla\eta)
    \left(
      \frac{1}{\sqrt{1+\norma{w_1*\nabla\eta}^2}}
      -
      \frac{1}{\sqrt{1+\norma{w_2*\nabla\eta}^2}}
    \right).
  \end{eqnarray*}
  Using the elementary inequality $\modulo{(1+x^2)^{1/2} -
    (1+y^2)^{1/2}} \leq \modulo{x-y}$ we obtain:
  \begin{eqnarray}
    \nonumber
    & &
    v (w_1) - v (w_2)
    \\
    \nonumber
    & = &
    \kappa \, \frac{(w_1-w_2)*\nabla\eta}{\sqrt{1+\norma{w_1*\nabla\eta}^2}}
    +
    \kappa \, \frac{w_2*\nabla\eta}{\sqrt{1+\norma{w_2*\nabla\eta}^2}}
    \frac{\sqrt{1+\norma{w_2*\nabla\eta}^2}
      -
      \sqrt{1+\norma{w_1*\nabla\eta}^2}}
    {\sqrt{1+\norma{w_1*\nabla\eta}^2}}
    \\
    \nonumber
    & &
    \norma{v (w_1) - v (w_2)}_{\L\infty (\reali^n;\reali^n)}
    \\
    \nonumber
    & \leq &
    \kappa \, \norma{(w_1-w_2)*\nabla\eta}_{\L\infty (\reali^n;\reali^n)}
    +
    \kappa \,
    \modulo{
      \norma{w_1*\nabla\eta}_{\L\infty (\reali^n;\reali^n)}
      -
      \norma{w_2*\nabla\eta}_{\L\infty (\reali^n;\reali^n)}
    }
    \\
    \nonumber
    & \leq &
    2 \, \kappa \, \norma{(w_1-w_2)*\nabla\eta}_{\L\infty (\reali^n;\reali^n)}
    \\
    \label{eq:diffv}
    & \leq &
    2 \, \kappa \, \norma{w_1-w_2}_{\L1 (\reali^n;\reali)}
    \, \norma{\nabla\eta}_{\L\infty (\reali^n;\reali^n)} \, .
  \end{eqnarray}
  Compute the divergence of $v (w)$ as follows:
  \begin{equation}
    \label{eq:div_v}
    \div v (w)
    =
    \kappa
    \frac{w * (\Delta \eta)}{\sqrt{1 + \norma{w*\nabla \eta}^2}}
    \left(
      1
      -
      \frac{(w * \nabla \eta) \, (w * \nabla \eta)}{1 + \norma{w*\nabla \eta}^2}
    \right)
    =
    \kappa
    \frac{w * (\Delta \eta)}{\left(1 + \norma{w*\nabla \eta}^2\right)^{3/2}}.
  \end{equation}
  We now compute the gradient of~\eqref{eq:div_v}:
  \begin{displaymath}
    \nabla \left(\div v (w)\right)
    =
    \kappa
    \frac{w* \nabla\Delta\eta}{\left(1+\norma{w*\nabla\eta}^2\right)^{3/2}}
    -
    3\kappa\,
    (w*\Delta\eta)
    \frac{w * \nabla^2\eta}{\left(1+\norma{w*\nabla\eta}^2\right)^2}
    \frac{w * \nabla\eta}{\sqrt{1 + \norma{w * \nabla \eta }^2}}
  \end{displaymath}
  so that
  \begin{eqnarray}
    \nonumber
    \!\!\!\!\!
    & \! &
    \!\!\!\!
    \norma{\nabla \left(\div v (w)\right)}_{\L1 (\reali^n;\ \reali^n)}
    \\
    \nonumber
    \!\!\!\!\!
    & \!\leq &
    \!\!\!\!
    \kappa \norma{w * \nabla\Delta\eta}_{\L1 (\reali^n; \reali^n)}
    +
    3\kappa\,
    \norma{w*\Delta\eta}_{\L1 (\reali^n; \reali)}  \,
    \norma{w*\nabla^2\eta}_{\L\infty (\reali; \reali^{n\times n})}
    \\
    \label{eq:graddivv}
    \!\!\!\!\!
    & \!\leq &
    \!\!\!\!
    \kappa
    \norma{w}_{\L1 (\reali^n; \reali)} \!\!
    \left(
      \norma{\nabla\Delta\eta}_{\L1 (\reali^n; \reali^n)} \!
      +
      3
      \norma{\Delta\eta}_{\L1 (\reali^n;\reali)} \,
      \norma{\nabla^2\eta}_{\L\infty (\reali^n;\reali^{n\times n})}
      \norma{w}_{\L1 (\reali^n;\reali)} \!
    \right)\!.
  \end{eqnarray}
  Consider now
  \begin{eqnarray*}
    & &
    \div\left(v (w_1) - v (w_2)\right)
    \\
    & = &\!
    \kappa\,
    \frac{(w_1-w_2)*\Delta\eta}{\left(1+\norma{w_1*\nabla\eta}^2\right)^{3/2}}
    +
    \kappa(w_2*\Delta\eta)\!
    \left(
      \frac{1}{\left(1+\norma{w_1*\nabla\eta}^2\right)^{3/2}}
      -
      \frac{1}{\left(1+\norma{w_2*\nabla\eta}^2\right)^{3/2}}
    \right).
  \end{eqnarray*}
  Using the inequality $\modulo{(1+x^2)^{-3/2} - (1+y^2)^{-3/2}} \leq
  \frac{48}{25 \sqrt{5}} \modulo{x-y}$,
  \begin{eqnarray}
    \nonumber
    & &
    \norma{\div\left(v (w_1) - v (w_2)\right)}_{\L1 (\reali^n; \reali)}
    \\
    \nonumber
    & \leq &
    \kappa\,
    \norma{(w_1-w_2)*\Delta\eta}_{\L1 (\reali^n; \reali)}
    +
    \frac{48}{25\sqrt{5}} \, \kappa \,
    \norma{w_2*\Delta\eta}_{\L\infty (\reali^n; \reali)}
    \norma{(w_1-w_2)*\nabla\eta}_{\L1 (\reali^n; \reali^n)}
    \\
    \nonumber
    & \leq &
    \kappa \, \norma{w_1-w_2}_{\L1 (\reali^n; \reali)} \,
    \norma{\Delta \eta}_{\L1 (\reali^n; \reali)}
    \\
    \nonumber
    & &
    +
    \frac{48}{25\sqrt{5}} \, \kappa \,
    \norma{w_2}_{\L\infty (\reali^n; \reali)}
    \norma{\Delta\eta}_{\L1 (\reali^n; \reali)}
    \norma{w_1-w_2}_{\L1 (\reali^n; \reali)}
    \norma{\nabla\eta}_{\L1 (\reali^n; \reali^n)}
    \\
    \label{eq:divdiffv}
    & \leq &
    \kappa \,
    \norma{w_1-w_2}_{\L1 (\reali^n; \reali)}
    \norma{\Delta\eta}_{\L1 (\reali^n; \reali)}
    \left(
      1 + \frac{48}{25\sqrt{5}} \,  \norma{\nabla\eta}_{\L1 (\reali^n; \reali^n)}
      \norma{w_2}_{\L\infty (\reali^n; \reali)}
    \right) \,.
  \end{eqnarray}
  Setting $K$ as in~\eqref{eq:4}, the inequalities above become:
  \begin{displaymath}
    \begin{array}{@{}rcll@{}}
      \norma{v (w)}_{\L\infty (\reali^n; \reali^n)}
      & \!\! \leq \! &
      K \norma{w}_{\L1 (\reali^n; \reali)}
      & \mbox{from~\eqref{eq:v}}
      \\
      \norma{\nabla v (w)}_{\L\infty (\reali^n; \reali)}
      & \!\! \leq \! &
      K \, \norma{w}_{\L\infty (\reali^n; \reali)}
      & \mbox{from~\eqref{eq:gradv}}
      \\
      \norma{v (w_1) - v (w_2) }_{\L\infty (\reali^n;\reali^n)}
      & \!\! \leq \! &
      K \,  \norma{w_1 -  w_2}_{\L1 (\reali^n; \reali)}
      & \mbox{from~\eqref{eq:diffv}}
      \\
      \norma{\nabla \left( \div v (w)\right)}_{\L1 (\reali^n; \reali^n)}
      & \!\! \leq \! &
      K \, \left( 1 + K \, \norma{w}_{\L1 (\reali^n; \reali)} \right)
      \norma{w}_{\L1 (\reali^n; \reali)}
      & \mbox{from~\eqref{eq:graddivv}}
      \\
      \norma{\div \left(v (w_1) - v (w_2)\right) }_{\L1
        (\reali^n;\reali)}
      & \!\! \leq \! &
      K \, \left( 1 + K \, \norma{w_2}_{\L\infty (\reali^n; \reali)} \right)
      \norma{w_1 - w_2}_{\L1 (\reali^n; \reali)}
      & \mbox{from~\eqref{eq:divdiffv}},
    \end{array}
  \end{displaymath}
  completing the proof.
\end{proofof}

\small{

  \bibliography{hyppara}

  \bibliographystyle{abbrv}

}

\end{document}